\documentclass[10pt]{article}
\usepackage[french]{babel}
\usepackage[latin1]{inputenc}
\usepackage[T1]{fontenc}
\usepackage{amssymb}
\usepackage{amsmath}
\usepackage{amsthm}
\usepackage{amscd}
\input xy
\xyoption{all}

\title{Foncteurs de division et structure de $I^{\otimes
    2}\otimes\Lambda^n$ dans la cat\'egorie $\mathcal{F}$}

\date{21 juillet 2006}

\author{Aur\'elien DJAMENT}

\newcommand{\C}{{\mathcal{C}}}

\newcommand{\F}{{\mathcal{F}}}
\newcommand{\E}{{\mathcal{E}}}

\newcommand{\FF}{{\mathbb{F}_2}}
\newcommand{\col}{{\rm colim}\,}

\newcommand{\go}{{\guillemotleft\,}}
\newcommand{\gf}{{\,\guillemotright\,}}

\newtheorem{thm-intro}{Théorème}

\newtheorem{theo}{Théorème}[section]
\newtheorem{pr}[theo]{Proposition}
\newtheorem{cor}[theo]{Corollaire}
\newtheorem{lm}[theo]{Lemme}
\newtheorem{thdef}[theo]{Théorème et définition}
\newtheorem{prdef}[theo]{Proposition et définition}

\theoremstyle{definition}
\newtheorem{defi}[theo]{Définition}
\newtheorem{nota}[theo]{Notation}

\newtheorem{conv}[theo]{Convention}

\theoremstyle{remark}
\newtheorem{rem}[theo]{Remarque}
\newtheorem{ex}[theo]{Exemple}

\begin{document}

\maketitle
 
%\section*{Introduction}

\begin{abstract} Nous démontrons que dans la catégorie $\F$ des foncteurs
entre espaces vectoriels sur $\FF$, le produit tensoriel entre le
second foncteur injectif standard non constant $V\mapsto{\mathbb{F}_2}^{(V^*)^{\oplus 2}}$ et un foncteur puissance
extérieure est artinien. Seul était antérieurement connu le caractère
artinien de cet injectif ; notre résultat constitue une étape pour l'étude
du troisième foncteur injectif standard non constant de $\F$.

Nous utilisons le foncteur de division par le foncteur identité et des
considérations issues de la théorie des représentations modulaires des
groupes symétriques pour obtenir ce théorème par la détection de
facteurs de composition convenables.

\begin{center}{\bf Abstract}
\end{center}

We prove that, in the category $\F$ of functors
  between $\FF$-vector spaces, the tensor product between the second
  non constant standard injective functor
  $V\mapsto{\mathbb{F}_2}^{(V^*)^{\oplus 2}}$  and an exterior power functor
  is artinian. The only case known to date was the artinian character
  of this injective ; our result is a step in the study of the third
  non constant standard injective of $\F$. 
 
We use the division functor by the identity functor and facts from modular representation theory of the symmetric groups to obtain
this theorem by detecting suitable composition factors.
\end{abstract}

\smallskip

\noindent
\begin{small}{\em Mots-clés : } cat\'egories de foncteurs, repr\'esentations
modulaires, foncteurs de division, filtration de Krull.

\noindent
{\em Classification math. : } 16P60, 18A25, 20B30, 20C20.
\end{small}

\section*{Introduction}

Cet article s'intéresse aux objets injectifs de la
catégorie $\F$ des foncteurs de la catégorie $\E^f$ vers la catégorie
$\E$, où l'on note $\E$ la catégorie des espaces vectoriels sur le
corps $\FF$ à deux éléments et $\E^f$ la sous-catégorie pleine de $\E$
dont les objets sont les espaces de dimension finie. Cette catégorie
s'est révélée fondamentale en algèbre et en topologie. Ainsi, Suslin a
démontré dans l'appendice de \cite{FFSS} que l'on peut calculer
certains groupes d'homologie stable des groupes linéaires sur $\FF$ à partir de
groupes d'extensions dans la catégorie $\F$, sur lesquels on dispose
de nombreux résultats (cf. \cite{FFSS} et \cite{Franjou}). La
catégorie $\F$ a été étudiée
systématiquement depuis les travaux de Henn, Lannes et Schwartz
montrant les liens étroits entre les foncteurs analytiques de $\F$ et
les modules instables sur l'algèbre de Steenrod (cf. \cite{HLS} et
\cite{LS}). L'image de la cohomologie d'un $2$-groupe abélien
élémentaire $V$ par le foncteur fondamental de \cite{HLS} se trouve
être l'objet injectif standard de $\F$ associé au $\FF$-espace
vectoriel $V$. Cela constitue l'une des motivations pour étudier ces objets injectifs.

Nous fournissons une étape supplémentaire dans l'étude de la {\em
  conjecture artinienne} (discutée en détails dans \cite{GP-gal} et \cite{these}) selon laquelle le foncteur injectif
$I^{\otimes d}$ est un objet artinien pour tout entier naturel $d$, où $I$ désigne le premier objet
injectif standard non constant de $\F$
(cf. §\,\ref{par-gll}). Celle-ci n'est démontrée que pour $d\leq
2$. Elle est triviale pour $d=0$, facile pour $d=1$ et beaucoup plus
difficile pour $d=2$. De fait, bien que la catégorie $\F$ possède de nombreuses propriétés de
régularité et que ses objets simples soient connus --- ils
sont paramétrisés par les représentations simples sur $\FF$ des différents
groupes symétriques, l'étude de ses objets de longueur infinie s'avère
d'une grande complexité, liée au problème de la compréhension {\em
  globale} de ces représentations. Les cas particuliers connus
ont requis l'utilisation d'outils puissants, qui sont pourtant encore insuffisants pour une approche générale.

Afin d'aborder l'étude de $I^{\otimes d}$, on commence par examiner les
sous-foncteurs obtenus en tensorisant par $I^{\otimes d-1}$ la
filtration polynomiale de $I$, dont les quotients font apparaître les
puissances extérieures.  Ainsi, pour établir le caractère artinien
de $I^{\otimes 2}$ dans \cite{P-art2}, Powell s'est appuyé sur le
résultat, dû à Piriou (\cite{Piriou}), selon lequel les foncteurs $I\otimes\Lambda^n$
sont artiniens. Nos résultats sur la structure des foncteurs $I^{\otimes
  2}\otimes\Lambda^n$, qui aboutiront au théorème suivant, constituent
donc une étape pour  l'étude de $I^{\otimes 3}$. 

\begin{thm-intro}\label{thp} Pour tout entier naturel $n$, le foncteur
  $I^{\otimes 2}\otimes\Lambda^n$ est artinien.
\end{thm-intro}

L'outil principal, et nouveau pour ce type de problème, que nous employons, est l'endofoncteur
$(-:\Lambda^1)$ de $\F$
adjoint à gauche à $-\otimes\Lambda^1$, appelé {\em division par
  $\Lambda^1$}, qui suggère une approche par récurrence du théorème
\ref{thp} à partir des isomorphismes naturels $(I^{\otimes
  2}\otimes\Lambda^n :\Lambda^1)\simeq I^{\otimes
  2}\otimes\Lambda^{n-1}$ et du cas $n=0$, traité par Powell.

Par des considérations issues de la théorie des représentations modulaires des
groupes symétriques,
nous donnons des renseignements sur l'effet de ce foncteur sur certains
objets simples de la catégorie $\F$. Nous montrons ensuite comment en
déduire des renseignements sur les facteurs de
composition de foncteurs analytiques dont on contrôle la division par $\Lambda^1$. Nous les appliquons à
des sous-foncteurs de $I^{\otimes 2}\otimes\Lambda^n$ pour obtenir le
pas de la récurrence. 

Précisément, nous introduisons des
sous-foncteurs explicites $L^2_n$ et $D^2_n$ de $I^{\otimes
  2}\otimes\Lambda^n$ qui permettent de réduire le théorème \ref{thp}
à l'énoncé suivant.

\begin{thm-intro}\label{thld} Pour tout entier naturel $n$, les
  foncteurs $L^2_n$ et $D^2_n$ sont artiniens.
\end{thm-intro}

En fait, nous donnons une estimation de la \go taille\gf des foncteurs
$L^2_n$, $D^2_n$ et  $I^{\otimes 2}\otimes\Lambda^n$ beaucoup plus
précise que leur caractère artinien ; elle est étroitement liée à la
filtration de Krull de $\F$ --- cf. \cite{these}.

Nous obtiendrons le théorème \ref{thld} à partir d'une proposition (\ref{pcruc}) affirmant qu'un
sous-foncteur $X$ de $L^2_n$  ayant la propriété que
l'image du morphisme $(X:\Lambda^1)\to (L^2_n:\Lambda^1)$ induit par l'inclusion est
suffisamment grosse est égal à $L^2_n$. Le pas de la récurrence s'en
déduit par un argument formel, car le foncteur $(L^2_n:\Lambda^1)$ est étroitement lié à $L^2_{n-1}$ ; le cas de
$D^2_n$ se traite pareillement.

\smallskip

Ce travail s'organise comme suit. La première section rappelle les
premières propriétés de la catégorie $\F$ et la construction de ses
objets simples. La deuxième
section est consacrée à l'étude des propriétés fondamentales de
l'endofoncteur $(-:\Lambda^1)$ de $\F$,  comparées à celles du foncteur
différence $\Delta$, l'un des outils de base les plus importants de la
catégorie $\F$. Le foncteur $(-:\Lambda^1)$ est un quotient du
foncteur différence ; contrairement à celui-ci, il n'est pas exact,
mais il possède d'autres propriétés qui le rendent bien plus
maniable que $\Delta$, souvent  trop \go gros\gf pour permettre des
calculs raisonnables, même sur les foncteurs assez élémentaires. La troisième section introduit une classe de foncteurs simples sur laquelle le
foncteur de division par $\Lambda^1$ se trouve assez maniable. La quatrième montre comment détecter la présence de facteurs de
composition dans un sous-foncteur $G$ d'un foncteur $F$ à partir de
celle de certains facteurs dans l'image de $(G:\Lambda^1)\to
(F:\Lambda^1)$. Dans la dernière section, après quelques préliminaires
sur la structure du foncteur
$I^{\otimes 2}$, nous construisons explicitement les foncteurs $L^2_n$ et $D^2_n$ qui permettent de
filtrer les $I^{\otimes 2}\otimes\Lambda^n$. Nous établissons les
propriétés de leurs facteurs de
composition et de leur division par $\Lambda^1$ grâce auxquelles on peut mener à bien l'argument de
récurrence des théorèmes \ref{thp} et \ref{thld}.

\smallskip

Cet article expose une partie des résultats obtenus par l'auteur dans
son doctorat \cite{these}, qui contient une
généralisation du théorème \ref{thp} au cas de $I^{\otimes 2}\otimes
F$, où $F$ est un foncteur fini. Celle-ci s'effectue au prix de
constructions plus techniques liées aux représentations des groupes
symétriques, et de l'utilisation occasionnelle d'autres outils, également introduits
dans \cite{these}, qui nous ont
permis d'obtenir aussi cette généralisation par une méthode n'utilisant pas
la division par $\Lambda^1$ (cf. \cite{art2} et \cite{art3}). Les deux
approches sont reliées par la nécessité, pour progresser
dans l'étude de la conjecture artinienne, d'une compréhension fine des
{\em foncteurs co-Weyl} (définis dans \cite{Po1}) et de phénomènes
\go globaux\gf dans les représentations modulaires des groupes
symétriques ou linéaires. 

\paragraph*{Notations et conventions}
\begin{itemize}\item Dans la suite, $\FF$ étant le seul corps de base
  que nous considérerons, nous nommerons simplement {\em espace
    vectoriel} un espace vectoriel sur $\FF$. 
\item Soit $E$ un ensemble. Nous noterons $\FF[E]$
  l'espace vectoriel  somme directe de copies de $\FF$ indexées par $E$. On peut voir
  l'assocation $E\mapsto\FF[E]$ comme un foncteur de la catégorie des ensembles vers $\E$.

Nous noterons
  $[e]$ l'élément de la base canonique de $\FF[E]$
  associé à un élément $e$ de $E$.
\item Si $G$ est un groupe, $\FF[G]$ désignera l'algèbre de $G$ sur $\FF$.
\item Nous noterons $V^*$ le dual d'un espace vectoriel $V$.
\item Nous désignerons par $\mathbf{Mod}_A$ la catégorie des modules à
  droite sur un anneau $A$.
\item Nous noterons ${\rm Ob}\,\C$ la classe des objets d'une
  catégorie $\C$. L'ensemble des morphismes d'un objet $A$ vers un objet $B$ de $\C$ sera noté ${\rm hom}_\C (A,B)$, ou
  ${\rm hom}(A,B)$ si nulle ambiguïté n'est possible.

Enfin, $\C^{op}$ désignera la catégorie opposée de $\C$.
\item La sous-catégorie pleine de $\F$ des foncteurs prenant des
  valeurs de dimension finie sera notée $\F^{df}$.
\item Nous désignerons par $\Lambda^n$ (resp. $T^n$, $S^n$, $\Gamma^n$) le foncteur
  $n$-ième puissance extérieure (resp. tensorielle, symétrique,
  divisée).
\item On note $\mathbb{Z}$ l'ensemble des entiers relatifs,
  $\mathbb{N}$ l'ensemble des entiers naturels et $\mathbb{N}^*$
  l'ensemble des entiers strictement positifs.
\end{itemize}

\tableofcontents

\section{Préliminaires}\label{sect-p}

\subsection{Rappels généraux sur la catégorie $\F$}\label{par-gll}

  Pour les définitions et résultats que nous rappelons dans ce paragraphe, et
dont nous ferons un usage fréquent, nous renvoyons le lecteur à \cite{K1}, \cite{Po2}
et \cite{LS} par exemple. Pour ce qui concerne les résultats de base
sur les catégories abéliennes que nous utilisons (souvent
implicitement), on peut se référer à \cite{Gab}.

\begin{prdef}\label{strglf}\begin{enumerate}\item La catégorie $\F$
    est abélienne, elle possède des limites et des colimites, qui se
    calculent au but ; les
    colimites filtrantes y sont exactes. De plus, le produit tensoriel
    de $\E$ induit une structure monoïdale symétrique sur  $\otimes :
    \F\times\F\to\F$, le {\em produit tensoriel}.
\item On définit un foncteur $(\E^f)^{op}\to\mathcal{F}\quad
  V\mapsto P_V$ par la formule $P_{V}(E)=\mathbb{F}_{2}[{\rm hom}_{\E_{f}}(V,E)]$. On note également $P=P_{\mathbb{F}_2}$. 
\item On a un isomorphisme ${\rm hom}_\mathcal{F}(P_V,F)\simeq
  F(V)$ naturel en $V\in {\rm Ob}\,\E^f$ et en $F\in\mathcal{F}$. En particulier,
  les foncteurs $P_V$ sont des objets projectifs de $\mathcal{F}$, appelés {\em foncteurs projectifs standard}
de $\F$.
\item Les foncteurs $P_{\mathbb{F}_2^{\oplus n}}$ ($n\in\mathbb{N}$) forment une famille de générateurs projectifs de $\mathcal{F}$.
\item On a un isomorphisme $P_V\otimes P_W\simeq P_{V\oplus W}$
  naturel en $V, W\in {\rm Ob}\,\E^f$. En particulier, $P_V\simeq
  P^{\otimes\dim V}$.
\end{enumerate}
\end{prdef}

\begin{prdef}\label{dualinj}\begin{enumerate}\item On définit un
    foncteur \\$D : \mathcal{F}^{op}\to\mathcal{F}$, appelé {\em dualité}, par
    $DF(V)=F(V^*)^*$.
\item Le foncteur $D$ induit une
  équivalence de catégories entre $(\F^{df})^{op}$ et $\F^{df}$.
\item Le foncteur $D$ est exact et fidèle ; il commute au produit tensoriel.
\item On a un isomorphisme ${\rm hom}(F,DG)\xrightarrow{dual_{F,G}}
  {\rm hom}(G,DF)$ naturel en les foncteurs $F$ et $G$.
\end{enumerate}
\end{prdef}

\begin{defi}
\label{adorth}
Un foncteur $F\in {\rm Ob}\,\mathcal{F}$ est dit {\em auto-dual} s'il est muni d'un isomorphisme $F\xrightarrow{d} DF$ tel que $dual_{F,F}(d)=d$.

Dans un tel foncteur, on définit l'{\em orthogonal} d'un sous-objet $A$ par
$$A^\perp=im\,\big(D(X/A)\xrightarrow{D\pi}DX\xrightarrow{d^{-1}} X\big),$$
où $\pi$ désigne la projection $X\twoheadrightarrow X/A$.
\end{defi}

\begin{rem} 
\begin{enumerate}
\item Un foncteur auto-dual est à valeurs de dimension finie.
\item Le produit tensoriel de deux foncteurs auto-duaux est auto-dual.
\end{enumerate}
\end{rem} 

\begin{ex}\label{exadl} Les foncteurs $\Lambda^i$ sont auto-duaux. Par
  conséquent, tout produit tensoriel de puissances extérieures est auto-dual.
\end{ex}

Il est usuel, en raison des liens avec les modules instables sur
l'algèbre de Steenrod (cf. \cite{HLS} et \cite{LS}), de s'intéresser plutôt aux objets injectifs de
$\F$ qu'à ses objets projectifs ; la dualité définie précédemment
permet de passer d'un point de vue à l'autre. 

\begin{prdef}\begin{enumerate}\item Pour $V\in {\rm Ob}\,\E^f$, on pose $I_V=DP_{V^*}$. On définit aussi
  $I=I_{\mathbb{F}_2}=DP$.
\item On a un isomorphisme ${\rm hom}_\mathcal{F}(F,I_V)\simeq
  F(V)^*$ naturel en $V\in {\rm Ob}\,\E^f$ et en $F\in {\rm Ob}\,\mathcal{F}$. En
  particulier, les foncteurs $I_V$ sont des objets injectifs de $\mathcal{F}$ ; on les appelle
  {\em foncteurs injectifs standard} de $\F$.
\item Les foncteurs $I_{\mathbb{F}_2^{\oplus n}}$ ($n\in\mathbb{N}$) forment une famille de cogénérateurs injectifs de $\mathcal{F}$.
\item On a un isomorphisme $I_V\otimes I_W\simeq I_{V\oplus W}$
  naturel en $V, W\in {\rm Ob}\,\E^f$.  En particulier, $I_V\simeq
  I^{\otimes\dim V}$.
\end{enumerate}
\end{prdef}

\begin{defi}[Foncteurs de décalage]
Pour tout $V\in {\rm Ob}\,\E^f$, on définit un endofoncteur $\Delta_V$ de $\mathcal{F}$ par $\Delta_V(F)=F\circ (V\oplus -)$. Noter que cela définit même un foncteur de $\E^f$ vers la catégorie des endofoncteurs de $\mathcal{F}$.
\end{defi}

On constate que ces foncteurs sont exacts et qu'ils commutent aux
limites, aux colimites et au produit tensoriel. De plus, on a des
isomorphismes naturels $\Delta_V\circ\Delta_W\simeq\Delta_{V\oplus
  W}$. 

\begin{pr}\label{prad1} Il existe un isomorphisme 
\begin{equation}
\label{adjel}
{\rm hom}_{\mathcal{F}}(A,\Delta_V B)\simeq {\rm
  hom}_{\mathcal{F}}(A\otimes P_V,B)
\end{equation}
naturel en $A,B\in {\rm Ob}\,\mathcal{F}$ et en $V\in {\rm Ob}\,\E^f$. En particulier, pour tout espace vectoriel $V$ de dimension finie, le
foncteur $\Delta_V$ est adjoint à droite au foncteur $-\otimes P_V$.
\end{pr}

On en déduit par dualité :

\begin{cor}\label{prad2} Il existe un isomorphisme 
\begin{equation}
\label{adjel2}
{\rm hom}_{\mathcal{F}}(\Delta_V  A, B)\simeq {\rm
  hom}_{\mathcal{F}}(A, I_V\otimes B)
\end{equation}
naturel en $A,B\in {\rm Ob}\,\mathcal{F}$ et en $V\in {\rm Ob}\,\E^f$. En particulier, pour tout espace vectoriel $V$ de dimension finie, le
foncteur $\Delta_V$ est adjoint à gauche au foncteur $-\otimes I_V$.
\end{cor}

\begin{pr}
Notons $i : \E\to\mathcal{F}$ le foncteur, exact et fidèle, qui à un
espace vectoriel $V$ associe le foncteur constant en $V$ (dans la
suite, nous noterons simplement $V$ pour $i V$). Le foncteur
$\mathcal{F}\to\E\quad F\mapsto F(0)$ d'évaluation en $0$ est adjoint à gauche et à droite
à $i$.

Le foncteur constant $i F(0)$ est naturellement facteur direct du
foncteur $F$ de $\F$.  
\end{pr}

\begin{nota}\label{nbar} Soit $F\in
{\rm Ob}\,\mathcal{F}$. Nous noterons $\bar{F}$ le conoyau de l'inclusion canonique
$F(0)\hookrightarrow F$, de sorte qu'on a un scindement naturel $F\simeq F(0)\oplus\bar{F}$.
\end{nota}

Nous utiliserons surtout cette notation pour le projectif $\bar{P}$ et
l'injectif $\bar{I}$.

\begin{cor}
Notons $\Delta : \mathcal{F}\to\mathcal{F}$ le foncteur conoyau de
l'inclusion canonique $id_\mathcal{F}\to\Delta_{\mathbb{F}_2}$, appelé
{\em foncteur différence} de $\F$. On a un scindement naturel
$\Delta_{\mathbb{F}_2}\simeq id\oplus\Delta$. De plus, le foncteur
$\Delta$ commute à la dualité ; il est adjoint à droite à $-\otimes\bar{P}$ et à gauche à $-\otimes\bar{I}$.
\end{cor}

\begin{defi}
\begin{enumerate}
\item Soit $n\in\mathbb{Z}\cup\{-\infty\}$. On note $\mathcal{F}_n$ la
  sous-catégorie pleine de $\mathcal{F}$ formée des foncteurs $F$ tels
  que $\Delta^{n+1}F=0$ si $n\geq 0$, réduite à $\{0\}$ sinon.
\item Un foncteur $F$ est dit {\em polynomial} s'il existe $n$ tel que $F$ appartient à $\mathcal{F}_n$. Le plus petit $n$ ayant cette propriété s'appelle le {\em degré} de $F$, on le note $\deg\,F$.
\item Un foncteur est dit {\em analytique} s'il est réunion de sous-foncteurs polynomiaux. On désigne par $\mathcal{F}_\omega$ la sous-catégorie pleine de $\mathcal{F}$ formée des foncteurs analytiques.
\end{enumerate}
\end{defi}

Dans la suite, nous nommerons simplement {\em finis} les objets de
longueur finie.

\begin{prdef}\label{POL}
\begin{enumerate}
\item Un foncteur est fini si et seulement s'il est polynomial et à valeurs de dimension finie.
\item Le dual d'un foncteur polynomial est polynomial de même degré.
\item Un foncteur est localement fini si et seulement s'il est analytique.
\item Les foncteurs $I_V$ sont analytiques.%\,\footnote{Ce n'est pas en revanche le cas de $P_V$ si $\dim\,V>0$.}.
\item Soit $n\in\mathbb{Z}\cup\{-\infty\}$. La sous-catégorie
  $\mathcal{F}_n$ de $\mathcal{F}$ est épaisse. De plus, l'inclusion
  $\mathcal{F}_{n}\xrightarrow{i_n}\mathcal{F}$ possède un adjoint à
  droite (resp. à gauche) noté $p_n$ (resp. $q_n$). Par abus, on
  désignera encore par $p_n$ (resp. $q_n$) l'endofoncteur $i_n p_n$
  (resp. $i_n q_n$) de $\mathcal{F}$.  
\item Soit $F$ un objet de $\F$. Via la coünité (resp. l'unité) de
  l'adjonction, on peut voir naturellement $p_n(F)$ comme un
  sous-objet (resp. un quotient) de $F$ ; de plus la suite
  $(p_n(F))_n$ de sous-objets de $F$ est croissante. Sa réunion est
  $F$ si et seulement si $F$ est analytique ; en général, cette
  réunion est le plus grand sous-foncteur analytique de $F$. 
\end{enumerate}
\end{prdef}

\begin{pr}\label{pol2}\begin{enumerate}
\item Si $A$ et $B$ sont deux foncteurs polynomiaux, alors $A\otimes B$ est également polynomial. De plus,  $\deg\,(A\otimes B)=\deg\,A+\deg\,B$. 
\item Soient $n\in\mathbb{Z}$ et $A, B\in {\rm Ob}\,\mathcal{F}$. On a un isomorphisme naturel
$$p_n(A\otimes B)\simeq\sum_{i+j=n}p_i(A)\otimes p_j(B)\,.$$
\item\label{pdh} Pour $n\in\mathbb{Z}$, on définit un endofoncteur $p^{hom}_n$ de $\mathcal{F}$ par
\begin{equation}
\label{defhom}
p^{hom}_n=p_n/p_{n-1}\,.
\end{equation}
Il existe un isomorphisme naturel
$$p^{hom}_n(A\otimes B)\simeq \bigoplus_{i+j=n}p^{hom}_i(A)\otimes p^{hom}_j(B)\,.$$
\end{enumerate}
\end{pr}

Pour une démonstration, on pourra se reporter à \cite{Pir-these}.

\begin{cor} Le produit tensoriel de deux foncteurs analytiques est analytique.
\end{cor}

%On remarquera que les foncteurs $p_n$ et $q_n$ sont duaux.

\begin{defi}
Un foncteur $F$ est dit {\em homogène de degré $n$} (resp. {\em cohomogène de degré $n$}) s'il est polynomial de degré $n$ et si $p_{n-1}(F)=0$ (resp. $q_{n-1}(F)=0$). 
\end{defi}

Il revient au même de dire que $F$ est de degré $n$ de même que tous ses sous-objets (resp. tous ses quotients) non nuls. 

\begin{rem} Le foncteur $p^{hom}_n(F)$ est homogène de degré $n$ s'il
  est non nul.
\end{rem}

\begin{cor}\label{crch}
\begin{enumerate}
\item Tout sous-foncteur (resp. tout quotient) d'un foncteur homogène (resp. cohomogène) est homogène (resp. cohomogène).
\item Un produit tensoriel de foncteurs homogènes (resp. cohomogènes) est également homogène (resp. cohomogène).
\end{enumerate}
\end{cor}

\begin{defi}
Un foncteur est dit {\em de type fini} (resp. {\em de co-type fini}) --- en abrégé tf (resp. co-tf) --- lorsqu'il est quotient d'une somme directe finie de foncteurs $P_V$ (resp. qu'il se plonge dans une somme directe finie de foncteurs $I_V$). 
\end{defi}

Cette définition d'objet tf de $\F$ est équivalente à la notion
catégorique générale. Pour une présentation détaillée des différentes
notions de finitude utiles dans l'étude de la catégorie $\F$, nous
renvoyons à \cite{these}.

\begin{pr}
\begin{enumerate}
\item Tout foncteur tf (resp. co-tf) est à valeurs de dimension
  finie. 
\item Un foncteur est tf si et seulement si son dual est co-tf.
\item Un foncteur fini est tf et co-tf.
\item Un foncteur est co-tf si et seulement s'il est analytique et de socle fini.
\end{enumerate}
\end{pr}

\subsection{Foncteurs de Weyl et foncteurs simples}\label{sws}

L'action par permutation des facteurs du groupe symétrique $\Sigma_n$ sur le foncteur $n$-ième puissance
tensorielle $T^n$ permet de décrire les objets simples de la catégorie
$\F$ à partir de ceux des catégories $\mathbf{Mod}_{\mathbb{F}_{2}[\Sigma_{n}]}$,
$n$ parcourant $\mathbb{N}$. Ce paragraphe rappelle la description
explicite des simples de $\F$ obtenue par cette approche, qui passe par
l'intermédiaire des foncteurs de Weyl, plus maniables que les
foncteurs simples eux-mêmes.

\begin{nota}\label{notrs} Soit $n\in\mathbb{N}$. 
\begin{itemize}\item Nous noterons $s_n :
\F_n\to\mathbf{Mod}_{\mathbb{F}_{2}[\Sigma_{n}]}$ le foncteur ${\rm
  hom}_\F(T^n,.)$, qu'on munit d'une structure de
$\mathbb{F}_{2}[\Sigma_{n}]$-module à droite en faisant
agir $\Sigma_n$ à {\em gauche} sur $T^n$.
\item Nous noterons $r_n :
  \mathbf{Mod}_{\mathbb{F}_{2}[\Sigma_{n}]}\to\F_n$ le foncteur
  $-\underset{\mathbb{F}_{2}[\Sigma_{n}]}{\otimes}
  T^n$.
\end{itemize}
\end{nota}

\begin{pr}\label{prfrep} Soit $n\in\mathbb{N}$. 
\begin{enumerate}\item Le foncteur $r_n$ est adjoint à gauche à $s_n$.
\item Les foncteurs $r_n$ et $s_n$ induisent des équivalences de
  catégories réciproques l'une de l'autre
$$\F_{n}/\F_{n-1}\simeq\mathbf{Mod}_{\mathbb{F}_{2}[\Sigma_{n}]}.$$
\end{enumerate}
\end{pr}

Cette proposition est démontrée dans \cite{Pir-these}.

%\begin{rem}\label{rq-prrec} Si $A$ et $B$ sont deux foncteurs polynomiaux non nuls de degrés respectifs $n$ et $m$, $s_{n+m}(A\otimes B)$ s'identifie canoniquement au $\Sigma_{n+m}$-module induit (relativement à l'inclusion évidente $\Sigma_n\times\Sigma_m\hookrightarrow\Sigma_{n+m}$) par le $\Sigma_n\times\Sigma_m$-module produit tensoriel extérieur de $s_n(A)$ et $s_m (B)$. On peut voir cette propriété comme un succédané à la propriété exponentielle pour le foncteur gradué $(T^n)_n$.
%\end{rem}

\smallskip

Nous rappelons maintenant, dans le contexte de la
catégorie $\F$, les définitions et propriétés fondamentales de la
théorie des représentations des groupes symétriques.

\begin{defi}\label{defpart}
\begin{enumerate} 
\item Une {\em partition} est une suite décroissante $\lambda$
  d'entiers, indexée par $\mathbb{N}^*$, qui stationne en $0$.
\item La {\em longueur} d'une partition $\lambda$ est le plus grand entier $r$, noté $l(\lambda)$, tel que $\lambda_r>0$. Si $\lambda$ est identiquement nulle, on convient que $l(\lambda)=0$. Par la suite, on identifiera une partition $\lambda$ et le $n$-uplet $(\lambda_1,\dots,\lambda_n)$ si $n\geq l(\lambda)$.
\item Une partition $\lambda$ est dite {\em
  $2$-régulière} si
$\lambda_i>\lambda_{i+1}$ pour $1\leq i< l(\lambda)$ ; le corps de
base étant fixé à $\FF$, nous parlerons
par la suite simplement de partition {\em régulière}.
\item Le {\em degré} d'une partition $\lambda$ est l'entier positif $|\lambda|=\underset{i\in\mathbb{N}^*}{\sum}\lambda_i$. Une partition de $n\in\mathbb{N}$ est par définition une partition de degré $n$.
\item Soient $\lambda$ et $\mu$ deux partitions de même degré. On
  notera $\lambda\leq\mu$ si
$$\forall n\in\mathbb{N}^*\qquad\sum_{i=1}^n\lambda_i\leq\sum_{i=1}^n\mu_i\,.$$
\end{enumerate}
\end{defi}

\begin{nota}\label{notpl}
Si $\lambda=(\lambda_1,\dots,\lambda_r)$ est un $r$-uplet d'entiers,
on notera $\Lambda^{\lambda}$, ou
$\Lambda^{\lambda_1,\dots,\lambda_r}$, le foncteur 
$$\Lambda^{\lambda}=\Lambda^{\lambda_1,\dots,\lambda_r}=\underset{1\leq
  i\leq r}{\bigotimes}\Lambda^{\lambda_i}.$$
% Pour $1\leq i\leq r$, on définit $\lambda^-_i$ (resp. $\lambda^+_i$)\index{Nota}{$\lambda^+_i$, $\lambda^-_i$} par $(\lambda^-_i)_j=\lambda_j$ si $j\neq i$ et $(\lambda^-_i)_i=\lambda_i-1$ (resp. $(\lambda^+_i)_j=\lambda_j$ si $j\neq i$ et $(\lambda^+_i)_i=\lambda_i+1$).
%
%Pour alléger les notations, nous simplifierons des écritures du type $(\lambda^+_a)^-_b$ en $\lambda^{+,-}_{a,b}$.
%
%On étendra ces notations de façon évidente à une suite infinie $\lambda$ d'entiers stationnant en $0$.
\end{nota}

\begin{rem} La relation $\leq$ définit un ordre {\em partiel} sur les  partitions, appelé parfois {\em ordre de dominance}  (cf. \cite{James}, §\,3). 
\end{rem}

\begin{rem}%\begin{enumerate}
%\item Si $\lambda$ est une partition régulière de longueur $r$, $\lambda^+_i$ et $\lambda^-_i$ (cf notation \ref{notpl}) sont des partitions pour $1\leq i\leq r+1$ et $1\leq i\leq r$ respectivement.  
%\item 
Soient $\lambda$ et $\mu$ deux partitions de même degré. 
\begin{itemize}
\item L'assertion $\lambda\leq\mu$ est équivalente à
$$\forall n\in\mathbb{N}^*\qquad\sum_{i\geq n}\lambda_i\geq\sum_{i\geq n}\mu_i.$$
\item Si $\lambda\leq\mu$, alors $l(\lambda)\geq l(\mu)$.
%\item Si $\lambda\leq\mu$, alors $\lambda^+_a\leq\mu^+_a$ (resp. $\lambda^-_a\leq\mu^-_a$) pour tout entier $a$ pour lequel l'inégalité fait sens. 
\end{itemize}
%\end{enumerate}
\end{rem}

Soient $i$ et $j$ deux entiers naturels. Il existe un unique morphisme
non nul $\Lambda^i\otimes\Lambda^j\to\Lambda^{i+j}$, appelé {\em
  produit}, et un unique morphisme non nul
$\Lambda^{i+j}\to\Lambda^i\otimes\Lambda^j$, appelé {\em coproduit}. Ces morphismes sont duaux.

\begin{nota}\label{nfd}
 \begin{enumerate}
\item Soient $i$, $j$ et $t$ des entiers tels que $0\leq t\leq
  j$ et $i\geq 0$. On note $\theta_{i,j,t} :
  \Lambda^i\otimes\Lambda^j\to\Lambda^i\otimes\Lambda^t\otimes\Lambda^{j-t}\to\Lambda^{i+t}\otimes\Lambda^{j-t}$ la flèche composée du coproduit sur le deuxième facteur tensorisé par $\Lambda^i$ et du produit sur les deux premiers facteurs tensorisé par $\Lambda^{j-t}$. Par auto-dualité des puissances extérieures,  nous identifierons $D\theta_{i,j,t}$ à un morphisme 
$\Lambda^{i+t}\otimes\Lambda^{j-t}\to\Lambda^{i}\otimes\Lambda^{j}$.

Lorsqu'aucune confusion n'est possible, nous omettrons les indices pour les morphismes $\theta$ et $D\theta$.
\item Soit $\lambda$ une partition de longueur $r$. On note, pour $1\leq i\leq r-1$ et $ 1\leq t\leq\lambda_{i+1}$,
\begin{equation*}
%\label{dfpsi}
\psi^{i,t}_\lambda=\Lambda^{\lambda_1,\dots,\lambda_{i-1}}\otimes\theta_{\lambda_i,\lambda_{i+1},t}\otimes\Lambda^{\lambda_{i+2},\dots,\lambda_r}:\Lambda^\lambda\to\Lambda^{\lambda_1,\dots,\lambda_{i-1},\lambda_i+t,\lambda_{i+1}-t,\lambda_{i+2},\dots,\lambda_r}
\end{equation*}
puis $$\psi_\lambda=\underset{1\leq
  t\leq\lambda_{i+1}}{\bigoplus_{1\leq i\leq r-1}}\psi^{i,t}_\lambda :
\Lambda^\lambda\to\underset{1\leq t\leq\lambda_{i+1}}{\bigoplus_{1\leq
    i\leq
    r-1}}\Lambda^{\lambda_1,\dots,\lambda_{i-1},\lambda_i+t,\lambda_{i+1}-t,\lambda_{i+2},\dots,\lambda_r}.$$
\end{enumerate}
 \end{nota}

\begin{defi}
Soit $\lambda$ une partition. On définit le {\em foncteur de
  Weyl} associé à $\lambda$, noté
$W_\lambda$, %(ou encore $W_{\lambda_1,\dots,\lambda_n}$ pour $n\geq
            %l(\lambda)$)
par $$W_\lambda=ker\,\psi_\lambda\subset\Lambda^\lambda.$$
\end{defi}

\begin{rem} \label{remweyl} Si $i\geq j$,
$$W_{(i,j)}=\bigcap_{1\leq t\leq j}ker\,\theta_{i,j,t}$$
et pour une partition $\lambda$ de longueur quelconque $r$
$$W_{\lambda}=\bigcap_{i=1}^{r-1}\Lambda^{\lambda_1,\dots,\lambda_{i-1}}\otimes W_{(\lambda_i,\lambda_{i+1})}\otimes\Lambda^{\lambda_{i+2},\dots,\lambda_r}\subset\Lambda^{\lambda}.$$
\end{rem}

\begin{rem}
Pour toute partition $\lambda$, $W_\lambda$ est non nul, c'est donc un foncteur homogène de degré $|\lambda|$. En revanche, le foncteur défini par le noyau analogue à celui qui fournit $W_\lambda$ est nul sur une suite d'entiers qui n'est pas une partition (cf. \cite{PS}).
\end{rem}

Dans la suite, on désigne par {\em cosocle} d'un objet le quotient de
celui-ci par son {\em radical} (intersection des sous-objets stricts
maximaux). Au moins sur les objets finis, le cosocle est le plus grand
quotient semi-simple, c'est donc la notion duale du
{\em socle} (plus
grand sous-objet semi-simple).

\begin{thdef}[Objets simples de
  $\F$]\label{repf}
Soit $\lambda$ une partition \textbf{régulière}. 
\begin{enumerate}
\item \label{nrf} Le radical de $W_\lambda$ est donné par $${\rm
    rad}\,W_\lambda=W_\lambda\cap W_\lambda^\perp$$
(cf. définition \ref{adorth} et exemple \ref{exadl}).
%où $\Lambda^\lambda$ est muni de sa structure auto-duale canonique.
\item Le cosocle de $W_\lambda$ est un objet simple de $\mathcal{F}$, appelé
  {\em foncteur de Schur} associé
  à $\lambda$ et noté $S_\lambda$ pour $n\geq l(\lambda)$). Celui-ci est de degré $|\lambda|$ ; en particulier, $W_\lambda$ est cohomogène. 
\item Les foncteurs de Schur sont auto-duaux.
\end{enumerate}

De plus, les foncteurs de Schur associés aux partitions régulières forment un système complet  de représentants des objets simples de $\mathcal{F}$.
\end{thdef}

\begin{theo}\label{fclambda}
Les facteurs de composition de $\Lambda^\lambda$, où $\lambda$ est une
partition de longueur $r$ de $n\in\mathbb{N}$, sont :
\begin{enumerate}
\item les $S_\mu$, où $\mu$ parcourt les partitions régulières de $n$
 telles que $\mu\geq\lambda$,
\item des simples du type $S_\mu$ avec $|\mu|<n$, $l(\mu)<r$, $\mu_1\geq\lambda_1$ et $\mu_{r-1}\leq\lambda_r$.
\end{enumerate}

En outre, $S_\lambda$ est facteur de composition unique de $\Lambda^\lambda$.
\end{theo}

Nous renvoyons le lecteur à \cite{PS}, \cite{K2} et \cite{James} pour
les deux résultats précédents. Il pourra également consulter le
chapitre 3 de \cite{Pir-these} pour une exposition complète des
résultats fondamentaux sur les objets simples de $\F$.

\begin{nota}
Etant donnés une partition régulière $\lambda$ et un foncteur $F\in
{\rm Ob}\,\F$, nous abrégerons l'assertion {\em $S_\lambda$ est
  facteur de composition} (i.e. sous-quotient) {\em de $F$} en $\lambda\vdash F$.
\end{nota}

\section{La division par $\Lambda^1$ dans $\F$}\label{sdl}

Cette section expose les propriétés de base de l'endofoncteur
$(-:\Lambda^1)$, dit de {\em division par $\Lambda^1$}, que nous
utiliserons pour repérer certains \go bons\gf facteurs de compositions
dans des foncteurs convenables. Les principales vertus de ce foncteur
sont les suivantes :
\begin{itemize}\item il est exact à droite et possède un comportement
  agréable vis-à-vis du produit tensoriel (c'est une dérivation) ;
\item c'est un quotient du foncteur différence ; si $F$ est un objet
  fini, $(F:\Lambda^1)$ contient les facteurs de composition
  de degré maximal de $\Delta F$ et  supprime la plupart des facteurs de degré inférieur, rendant son
  calcul plus facile --- ainsi, nous verrons au paragraphe \ref{sdlw}
  que l'on peut calculer aisément la division par
  $\Lambda^1$ d'un foncteur de Weyl, contrairement à ce qui advient
  pour le foncteur différence ;
\item le foncteur de division
par $\Lambda^1$ généralise naturellement le foncteur de restriction
des $\Sigma_n$-modules vers les $\Sigma_{n-1}$-modules (cf. paragraphe
\ref{parlfs}), remarque qui rejoint la précédente ;
\item contrairement au foncteur différence, qui accroît la taille des
  foncteurs infinis, la division par $\Lambda^1$ diminue celle des
  foncteurs de co-type fini.
\end{itemize}

La principale difficulté occasionnée par l'emploi de ce foncteur,
comparé au foncteur différence, réside dans son inexactitude. Ainsi,
si $Y$ est un sous-quotient de $X$, $(Y : \Lambda^1)$ n'est pas
forcément un sous-quotient de $(X : \Lambda^1)$. Nous utiliserons
donc, dans la section \ref{sctdl}, des méthodes de détection de facteurs
de composition adaptées à des foncteurs seulement exacts à droite.

Signalons que les foncteurs de division ont été introduits par Lannes
dans \cite{Lannes} dans le cadre de la catégorie des modules instables sur l'algèbre
de Steenrod, intimement liée à $\F$. Dans \cite{Po2}, §\,$3$,
Powell met en évidence des liens étroits entre les foncteurs de division
considérés par Lannes et ceux de la catégorie $\F$.

%\subsection{Généralités sur les bifoncteurs $\mathbf{Hom}$ et $(-:-)$ de $\F$}
\subsection{Les bifoncteurs $\mathbf{Hom}$ et $(-:-)$
  de $\F$}

 \begin{prdef}\label{osf} Soit $F$ un objet de $\F$.
\begin{enumerate}\item L'endofoncteur $-\otimes F$ de $\F$ possède un
  adjoint à droite, noté \\
$\mathbf{Hom}(F,-)$, et appelé {\em
  foncteur hom interne de source $F$}. Ce foncteur est donc exact à gauche.
\item Si $F$ est objet est à
  valeurs de dimension finie, alors  $-\otimes F$ possède un adjoint à gauche, noté $(-:F)$ et
appelé {\em foncteur de division par $F$}. Ce foncteur est donc exact à droite.
\end{enumerate}
\end{prdef}

\begin{proof} Le foncteur $-\otimes F$ commute toujours aux colimites
  ; si $F$ est à valeurs de dimension finie, il commute également aux
  limites. La conclusion provient donc du théorème de Freyd.
\end{proof}

On obtient même ainsi des bifoncteurs $\mathbf{Hom} : \F^{op}\times\F\to\F$
et $(-:-) : \F\times (\F^{df})^{op}\to\F$. Ils sont liés par
l'isomorphisme naturel de dualité
\begin{equation}\label{eqdhid}\mathbf{Hom}(F,DG)\simeq D(G:DF).\end{equation}

L'isomorphisme (\ref{adjel}) de la proposition \ref{prad1} fournit
$\mathbf{Hom}(P_V,-)\simeq\Delta_V$, tandis que l'isomorphisme (\ref{adjel2})
du corollaire \ref{prad2} donne $(-:I_V)\simeq\Delta_V$.

On en déduit en particulier que les foncteurs hom internes et de
division commutent aux foncteurs différence et de décalage ; de plus,
si $F$ est un foncteur de type fini (resp. de co-type fini), alors
$\mathbf{Hom}(F,-)$ (resp. $(-:F)$) conserve $\F^{df}$ et les objets de tf
(resp. de co-tf), puisque $F$ est un quotient d'une somme
directe finie de $P_V$ (resp. un sous-objet d'une somme directe finie
de $I_V$).

\begin{pr}\label{homhi} Soient $X$ un objet de $\F$ et $A$ un objet homogène de
  degré $k$. Pour
  tout entier $n$, il existe un isomorphisme
  naturel $\mathbf{Hom}(A,p_n(X))\simeq p_{n-k}(\mathbf{Hom}(A,X))$.
\end{pr}

\begin{proof} C'est une conséquence formelle de l'adjonction
  entre $p_n$ et $i_n$ et de l'isomorphisme naturel 
$p_{n}(A\otimes B)\simeq A\otimes p_{n-k}(B)$ (cf. proposition/définition \ref{POL}).\end{proof}

\smallskip

Nous terminons ces généralités par une propriété relative aux {\em
  foncteurs exponentiels}, qui constituent un outil très commode pour mener à
bien des calculs sur des produits tensoriels dans $\F$
(cf. \cite{FFSS} et \cite{Franjou}).

\begin{defi}[Foncteurs exponentiels]  On appelle {\em foncteur
    exponentiel gradué} toute suite $(E^n)_{n\in\mathbb{N}}$ d'objets de $\F^{df}$ telle
  qu'il existe des isomorphismes
$$E^{n}(V\oplus W)\simeq\underset{i+j=n}{\bigoplus}E^{i}(V)\otimes
E^{j}(W)$$
naturels en les objets $V$ et $W$ de $\E^f$.
\end{defi}

\begin{ex} Les foncteurs $(\Lambda^n)_{n\in\mathbb{N}}$,
  $(S^n)_{n\in\mathbb{N}}$ et $(\Gamma^n)_{n\in\mathbb{N}}$ sont
  exponentiels gradués.
\end{ex}

\begin{pr}\label{prexpd}
Soient $(E^{n})_{n\in\mathbb{N}}$ un foncteur exponentiel gradué, $A$
et $B$ deux objets de $\F$. On a des isomorphismes naturels
$$\mathbf{Hom}(E^{n},A\otimes B)\simeq\bigoplus_{i+j=n}\mathbf{Hom}(E^{i},A)\otimes\mathbf{Hom}(E^{j},B)$$
et
$$(A\otimes B:E^{n})\simeq\bigoplus_{i+j=n}(A:E^{i})\otimes (B:E^{j})\,.$$
\end{pr}

\begin{proof}  Un argument de dualité permet de ne traiter que le
  premier cas. 

La structure exponentielle de $E$
  fournit, pour $i+j=n$, un morphisme $E^n\to E^i\otimes E^j$
  (coproduit), d'où un morphisme naturel 
$$\mathbf{Hom}(E^i,A)\otimes\mathbf{Hom}(E^j,A)\otimes E^n\to
(\mathbf{Hom}(E^i,A)\otimes E^i)\otimes (\mathbf{Hom}(E^j,A)\otimes
E^j)$$
$\to A\otimes B$,
où la dernière flèche est le produit tensoriel des deux morphismes
procurés par la coünité. Par adjonction, on en déduit un morphisme
naturel
$$\bigoplus_{i+j=n}\mathbf{Hom}(E^{i},A)\otimes\mathbf{Hom}(E^{j},B)\to\mathbf{Hom}(E^{n},A\otimes B)\,,$$
dont nous allons voir que c'est un isomorphisme. 

L'assertion analogue pour ${\rm hom}$ est démontrée, par
  exemple, dans \cite{FFSS}. Le cas général s'en déduit via les
  isomorphismes naturels
$$\mathbf{Hom}(F,G)(V)\simeq {\rm hom}(P_V,\mathbf{Hom}(F,G))\simeq
{\rm hom}(P_V\otimes F,G)\simeq {\rm hom}(F,\Delta_V G)$$
et la commutation des foncteurs $\Delta_V$ au produit tensoriel.
\end{proof}

\subsection{Propriétés générales des foncteurs $(-:\Lambda^1)$ et
  $\mathbf{Hom}(\Lambda^1,-)$}

Nous exposons dans cette section quelques propriétés générales des endofoncteurs $\mathbf{Hom}\,(\Lambda^1,-)$ et $(-:\Lambda^1)$ de
$\F$, qui jouent un rôle particulier parmi tous les foncteurs hom
internes et de division par un objet fini. Ces foncteurs sont {\em
  duaux} : il existe un isomorphisme
\begin{equation}\label{hidfleq} D(F :
\Lambda^1)\simeq\mathbf{Hom}\,(\Lambda^1,DF)\end{equation}
naturel en l'objet $F$
de $\F$ grâce à l'isomorphisme (\ref{eqdhid}), puisque le foncteur $\Lambda^1$ est auto-dual. 

\begin{pr}\label{deriv}
Les endofoncteurs $\mathbf{Hom}\,(\Lambda^1,-)$ et $(-:\Lambda^1)$ de
$\F$ sont des dérivations en ce sens qu'on a des isomorphismes
$$\mathbf{Hom}\,(\Lambda^1,F\otimes G)\simeq(\mathbf{Hom}\,(\Lambda^1,F)\otimes G)\oplus (F\otimes\mathbf{Hom}\,(\Lambda^1,G))$$
et
$$(F\otimes G:\Lambda^1)\simeq\big( (F:\Lambda^1)\otimes G\big)\oplus
\big(F\otimes (G:\Lambda^1)\big)$$ naturels en les objets $F$ et $G$
de $\F$.
\end{pr}

\begin{proof} Cela résulte de la proposition \ref{prexpd}.
\end{proof}

Dans la suite, nous noterons souvent $(F\otimes
G:\Lambda^1)\twoheadrightarrow (F:\Lambda^1)\otimes G$ sans plus de
précision la projection naturelle déduite de cette proposition. 

L'isomorphisme naturel $\Delta (F\otimes
G)\simeq (F\otimes\Delta G)\oplus (\Delta F\otimes G)\oplus (\Delta
F\otimes\Delta G)$ est à comparer à la proposition \ref{deriv}.

\begin{lm}
Soit $F$ un foncteur polynomial non nul de degré $d$. Il existe un morphisme
non nul de $F$ vers $T^d$.
\end{lm}

\begin{proof} Le foncteur $\Delta^d F$ est non nul, et constant parce que $\Delta^{d+1}
  F=0$, donc il existe un morphisme non nul de $F$ vers
  $\bar{I}^{\otimes d}$, par adjonction. Comme $\deg F\leq d$, ce morphisme est à valeurs dans
  le sous-objet $p_d(\bar{I}^{\otimes d})\simeq T^d$ de $\bar{I}^{\otimes d}$, d'où le lemme.\end{proof}

\begin{lm}\label{lmutii}
Si $F$ est un foncteur fini tel que $F(0)=0$ et $(F:\Lambda^1)=0$, alors $F=0$.
\end{lm}

\begin{proof} On a hom$\,(F,T^n)\simeq{\rm
    hom}\,((F:\Lambda^1),T^{n-1})=0$ pour tout $n\in\mathbb{N}^*$. On conclut par le lemme précédent.\end{proof}

\begin{pr}\label{prc-anhl}
Soit $X\in {\rm Ob}\,\mathcal{F}$. Le foncteur $(X:\Lambda^1)$ est nul
si et seulement si $\overline{\Delta^n X}$ est sans quotient fini non
nul pour tout $n\in\mathbb{N}$.
\end{pr}

\begin{proof} Supposons $(X:\Lambda^1)=0$. Le lemme
  précédent et l'exactitude à droite de
  $(-:\Lambda^1)$ impliquent que $\bar{X}$ n'a pas
  de quotient fini non nul. La commutation de
$(-:\Lambda^1)$ et $\Delta^n$ montre alors qu'il en est de même pour les $\overline{\Delta^n X}$.

Réciproquement, si les $\overline{\Delta^n X}$ sont sans quotient
fini non nul, il en est de même des $\overline{\Delta_{V}X}$ ($V\in
{\rm Ob}\,\E^{f}$), donc $(X:\Lambda^1)(V)^*\simeq {\rm
  hom}((X:\Lambda^1),I_V)\simeq {\rm hom}(X,\Lambda^1\otimes I_V)\simeq{\rm hom}(\Delta_{V}X,\Lambda^1)=0$,
ce qui achève la démonstration. \end{proof}

Cette proposition fournit les importants corollaires suivants.

\begin{cor}
\label{degfd}
Soit $F$ un foncteur fini non constant. Le foncteur $(F:\Lambda^1)$
est fini et non nul ; de plus, $\deg\, (F:\Lambda^1)=\deg F-1$.
\end{cor}

\begin{cor}\label{hil-ptf}
Pour tout foncteur injectif de co-type fini $X$ de $\F$, on a $(X:\Lambda^1)=0$.
\end{cor}

\begin{proof} Le foncteur $\bar{I}$ n'a pas de quotient fini non nul,
  comme il résulte par exemple du théorème $7.8$ de \cite{K2}, donc un foncteur injectif co-tf n'a pas de quotient fini non
  constant. Comme le foncteur différence $\Delta$ préserve les objets
  injectifs co-tf, la proposition \ref{prc-anhl} donne la conclusion.
\end{proof}

\begin{cor}\label{cr-comdl} Si $X$ est un objet injectif co-tf de
  $\F$, les endofoncteurs $-\otimes X$ et $(-:\Lambda^1)$ de $\F$
  commutent à isomorphisme canonique près.
\end{cor}

\begin{proof} Il s'agit d'une conséquence de la proposition
  \ref{deriv} et du corollaire \ref{hil-ptf}.
\end{proof}

La propriété suivante des foncteurs
$(-:\Lambda^1)$ et $\mathbf{Hom}\,(\Lambda^1,-)$  s'avère
fondamentale pour effectuer des calculs sur des foncteurs finis.

\begin{pr}
\label{fonddivl}
On a des suites exactes naturelles en $F\in {\rm Ob}\,\mathcal{F}$
\begin{equation}
\label{sef1}
0\to \mathbf{Hom}\,(\Lambda^1,F)\xrightarrow{\alpha_F}\Delta F\to \Delta^2 F  
\end{equation}
et
\begin{equation}
\label{sef2}
\Delta^2 F\to\Delta F\xrightarrow{\beta_F} (F:\Lambda^1)\to 0\,.
\end{equation}

De plus, si l'on note $\upsilon_{F}=\beta_F \alpha_F : \mathbf{Hom}\,(\Lambda^1,F)\to
(F:\Lambda^1)$, on a les résultats suivants.
\begin{enumerate}
\item \label{pfo1} Si $F$ est de degré $n$, alors
$$ker\,\beta_F\subset p_{n-2}(\Delta F)\quad\text{et}\quad ker\,\upsilon_{F}\subset p_{n-2}(\mathbf{Hom}\,(\Lambda^1,F))\simeq\mathbf{Hom}\,(\Lambda^1,p_{n-1}F).$$
\item \label{pfo2} Si $F$ est homogène, $\upsilon_{F}$ est injectif.
\item \label{pfo3} Si $F$ est cohomogène, $\upsilon_{F}$ est surjectif.
\item \label{pfo4} Si $F$ est homogène et cohomogène de degré $n$, $\mathbf{Hom}\,(\Lambda^1,F)$ et $(F:\Lambda^1)$, qui sont naturellement isomorphes via $\upsilon_{F}$, s'identifient
à $p_{n-1}^{hom}(\Delta F)$, qui est facteur direct de $\Delta F$, où la notation $p^{hom}_k$ est définie
dans la proposition \ref{pol2} par (\ref{defhom}).
\end{enumerate}
\end{pr}

\begin{proof}
On utilise la suite exacte usuelle
$$0\to\Lambda^1\to\bar{I}\to\bar{I}^{\otimes 2}$$
déduite de l'isomorphisme $p_1(\bar{I})\simeq\Lambda^1$ --- cf. \cite{K1}, lemme $4.12$
--- et sa duale 
$$\bar{P}^{\otimes 2}\to\bar{P}\to\Lambda^1\to 0$$
pour obtenir les suites
exactes (\ref{sef2}) et (\ref{sef1}). L'assertion (\ref{pfo1})
découle ensuite de la proposition \ref{homhi} et de ce que $\deg\Delta^2 F=\deg F-2$ si cet entier est
positif, $-\infty$ sinon ; (\ref{pfo2}) en résulte et
(\ref{pfo3}) se déduit de (\ref{pfo2}) par dualité.

Si $F$ est homogène et cohomogène, $\upsilon_F$ est un isomorphisme
par (\ref{pfo2}) et (\ref{pfo3}), d'où  le scindement $\Delta
F\simeq\mathbf{Hom}\,(\Lambda^1,F)\oplus ker\,\beta_F$, par définition de
$\upsilon_F$. Comme
$ker\,b_{F}\subset p_{n-2}(\Delta F)$ et $\mathbf{Hom}\,(\Lambda^1,F)\cap p_{n-2}(\Delta F)=0$, puisque $\mathbf{Hom}\,(\Lambda^1,F)$ est homogène, la composée $\mathbf{Hom}\,(\Lambda^1,F)
\hookrightarrow\Delta F\twoheadrightarrow p_{n-1}^{hom}(\Delta F)$ est un isomorphisme, ce qui achève de prouver (\ref{pfo4}).
\end{proof}

\begin{cor}
\label{cordiv}
Soit $\lambda=(\lambda_1,\dots,\lambda_r)$ une suite finie
d'entiers. Les foncteurs $\mathbf{Hom}\,(\Lambda^1,\Lambda^{\lambda})$ et $(\Lambda^{\lambda}:\Lambda^1)$ sont isomorphes à
$\underset{1\leq i\leq r}{\bigoplus}\Lambda^{\lambda^-_{i}}$.
\end{cor}

\begin{proof} Le corollaire \ref{deriv} montre qu'il suffit de vérifier
l'assertion pour $r=1$, auquel cas elle découle de la dernière
assertion de la proposition précédente, compte-tenu de $\Delta\Lambda^n=\Lambda^{n-1}$.\end{proof}

Nous terminons cette section par un
calcul explicite élémentaire.

\begin{ex} Nous allons déterminer
$\mathbf{Hom}\,(\Lambda^1,p_n\bar{I})$ et
$(p_n\bar{I}:\Lambda^1)$. Par le lemme $4.12$ de \cite{K1}, on a pour tout $n\in\mathbb{N}$ une suite exacte
\begin{equation*}
0\to p_{n-1} I\to p_{n}I\to\Lambda^{n}\to 0.
\end{equation*}

De plus, $\Lambda^k$ est le cosocle de $p_k\bar{I}$, de sorte que, par
la proposition précédente, $\upsilon_{p_k\bar{I}}$ est surjectif. On a
aussi des morphismes
$p_{n-1} I\otimes\Lambda^1\xrightarrow{a_n}p_n\bar{I}$, obtenus
en appliquant $p_n$ à l'unique morphisme non nul $I\otimes\Lambda^1\to\bar{I}$, tels que les diagrammes suivants commutent :
\begin{equation*}
\xymatrix{p_{n-1} I\otimes\Lambda^1\ar[d]_{a_n}\ar@{^{(}->}[r] &
  p_{n} I\otimes\Lambda^1\ar[d]^{a_{n+1}}\ar@{>>}[r] & \Lambda^n\otimes\Lambda^1\ar@{>>}[d] \\
p_{n}\bar{I}\ar@{^{(}->}[r] & p_{n+1}\,\bar{I}\ar@{>>}[r] & \Lambda^{n+1}
}\end{equation*}

On en déduit par adjonction des morphismes
$p_{n-1} I\xrightarrow{b_n}\mathbf{Hom}\,(\Lambda^1,p_{n}\bar{I})$
faisant commuter les diagrammes
\begin{equation*}
\xymatrix{p_{n-1} I\ar[d]_{b_n}\ar@{^{(}->}[r] & p_{n} I\ar[d]^{b_{n+1}}\ar@{>>}[r] & \Lambda^n\ar@{=}[d] \\
\mathbf{Hom}\,(\Lambda^1,p_{n}\bar{I})\ar@{^{(}->}[r] & \mathbf{Hom}\,(\Lambda^1,p_{n+1}\bar{I})\ar[r] & \mathbf{Hom}\,(\Lambda^1,\Lambda^{n+1}). 
}\end{equation*}

On note enfin $c_n$ le morphisme
$(p_n\bar{I}:\Lambda^1)\twoheadrightarrow
(\Lambda^n:\Lambda^1)=\Lambda^{n-1}$.

Nous allons montrer, par récurrence sur $n\in\mathbb{N}$, que les flèches $b_n$ et $c_n$ sont des
isomorphismes  et  le diagramme suivant commute.
\begin{equation*}
\xymatrix{\mathbf{Hom}\,(\Lambda^1,p_n\bar{I})\ar@{>>}[d]_{\upsilon_{p_{n}I}} &
  p_{n-1} I\ar[l]_-{\simeq}^-{b_n} \ar@{>>}[d] \\
(p_n\bar{I}:\Lambda^1)\ar[r]^(.55)\simeq_-{c_n} & \Lambda^{n-1}
}\end{equation*}

Pour $n=0$, l'assertion est évidente. Pour déduire l'assertion pour
$n+1$ de l'assertion pour $n$, on considère le diagramme commutatif
aux lignes exactes suivant.
\begin{equation*}
\xymatrix{p_{n-1}I\ar@{^{(}->}[r]\ar[d]_{b_n}^\simeq & p_n I
  \ar@{>>}[r]\ar[d]_{b_{n+1}} & \Lambda^n\ar@{=}[d] \\
\mathbf{Hom}\,(\Lambda^1,p_n\bar{I})\ar@{^{(}->}[r] \ar@{>>}[d]_{\upsilon_{p_n\bar{I}}}
  &
  \mathbf{Hom}\,(\Lambda^1,p_{n+1}\bar{I})\ar[r]\ar@{>>}[d]_{\upsilon_{p_{n+1}\bar{I}}}
  & \mathbf{Hom}\,(\Lambda^1,\Lambda^{n+1})\ar@{=}[d] \\
(p_n\bar{I}:\Lambda^1)\ar[d]_{c_n}^\simeq \ar[r] &
  (p_{n+1}\bar{I} :\Lambda^1)\ar@{>>}[r]^{c_{n+1}} & \Lambda^n \\
\Lambda^{n-1} & & 
}\end{equation*}

La commutation du carré en bas à droite entraîne que la flèche
horizontale centrale de droite est surjective, ce qui permet de
conclure quant à $\mathbf{Hom}\,(\Lambda^1,p_{n+1}\bar{I})$ en
appliquant le lemme des cinq à la partie supérieure du diagramme. Il
suffit donc d'établir la nullité de la flèche $(p_{n}\bar{I}:\Lambda^1)\to (p_{n+1}\bar{I}:\Lambda^1)$ induite par l'inclusion. Son image est de degré au plus $n-1$ et sans terme
constant (pour $n=1$, cela vient de ce que la suite exacte $0\to\Lambda^1\to p_{2}\bar{I}\to\Lambda^2\to 0$ est non scindée ; sinon cela découle de l'hypothèse de récurrence), mais c'est aussi un quotient de $\mathbf{Hom}\,(\Lambda^1,p_{n+1}\bar{I})\simeq p_{n}I$ ; du fait que
$p_{n}\bar{I}$ est cohomogène de degré $n$ (son cosocle étant $\Lambda^n$), cela entraîne la nullité de ladite image.
\end{ex}

\subsection{Liens formels avec les représentations des groupes
  symétriques}\label{parlfs}

Nous explicitons à présent en quoi les foncteurs $(-:\Lambda^1)$,
$\mathbf{Hom}(\Lambda^1,-)$ et $\Delta$ constituent 
des généralisations du foncteur de restriction pour les
représentations des groupes symétriques. Cela fait l'objet de la
proposition \ref{repdiv-c}. La proposition \ref{repdiv-p} montrera
qu'en un certain sens, le foncteur de division par $\Lambda^1$
constitue une meilleure généralisation que ses homologues de la proposition \ref{repdiv-c}.

\begin{conv} Dans ce paragraphe, $n$ désigne un entier naturel. On notera ${\rm Res} :
\mathbf{Mod}_{\FF[\Sigma_n]}\to\mathbf{Mod}_{\FF[\Sigma_{n-1}]}$ le
foncteur de restriction des scalaires et ${\rm Ind} : \mathbf{Mod}_{\FF[\Sigma_{n-1}]}\to\mathbf{Mod}_{\FF[\Sigma_n]}$ le foncteur d'induction.
\end{conv}

Les foncteurs $r_n$ et $s_n$ apparaissant ci-après sont ceux de la
notation \ref{notrs}.
%
%Les propositions suivantes précisent le lien entre division par $\Lambda^1$ et représentations des groupes symétriques ; comme la proposition  \ref{fonddivl} elles permettent en quelque sorte de pallier l'inexactitude de $(-:\Lambda^1)$.

\begin{pr}
\label{repdiv-c}
Pour tout entier $n$, le diagramme
\begin{equation}
\label{repdiv}
\xymatrix{\mathcal{F}_n\ar[d]_-{s_n}\ar[r]^-{(-:\Lambda^1)} & \mathcal{F}_{n-1}\ar[d]^{s_{n-1}}\\
\mathbf{Mod}_{\mathbb{F}_2[\Sigma_n]}\ar[r]^-{{\rm Res}} & \mathbf{Mod}_{\mathbb{F}_2[\Sigma_{n-1}]}
}
\end{equation}
commute à isomorphisme naturel près.

La même assertion vaut en remplaçant $(-:\Lambda^1)$ par $\mathbf{Hom}(\Lambda^1,-)$ ou $\Delta$.
\end{pr}

\begin{proof} On établit l'assertion relative à $\mathbf{Hom}(\Lambda^1,-)$. Les autres s'en déduisent par la proposition \ref{fonddivl}.

Grâce aux adjonctions de la proposition \ref{prfrep} et entre ${\rm Res}$ et
${\rm Ind}$, il suffit de montrer que le diagramme
$$\xymatrix{\mathcal{F}_n & \mathcal{F}_{n-1}\ar[l]_{-\otimes\Lambda^1}\\
\mathbf{Mod}_{\mathbb{F}_2[\Sigma_n]}\ar[u]^{r_n} & \mathbf{Mod}_{\mathbb{F}_2[\Sigma_{n-1}]}\ar[u]_{r_{n-1}}\ar[l]^-{{\rm Ind}}
}$$
commute à isomorphisme naturel près.

En effet, les propriétés d'associativité du produit tensoriel procurent
dans $\F$ un isomorphisme
$$(M\underset{\FF[\Sigma_{n-1}]}{\otimes}\FF[\Sigma_n])\underset{\FF[\Sigma_{n}]}{\otimes}T^n\simeq(M\underset{\FF[\Sigma_{n-1}]}{\otimes}
T^{n-1})\otimes\Lambda^1$$
naturel en l'objet $M$ de $\mathbf{Mod}_{\mathbb{F}_2[\Sigma_{n-1}]}$.
\end{proof}

\begin{pr}
\label{repdiv-p}
Pour tout entier $n$, le diagramme
\begin{equation}
\label{repdiv2}
\xymatrix{\mathbf{Mod}_{\mathbb{F}_2[\Sigma_n]}\ar[r]^(.45){{\rm Res}}\ar[d]_{r_n} & \mathbf{Mod}_{\mathbb{F}_2[\Sigma_{n-1}]}\ar[d]^{r_{n-1}} \\
\mathcal{F}_n\ar[r]^(.45){(-:\Lambda^1)} & \mathcal{F}_{n-1}
}
\end{equation}
commute à isomorphisme naturel près.
\end{pr}

\begin{proof} La tensorisation par $\Lambda^1$ procure un
  monomorphisme $\Sigma_{n-1}$-équivariant 
${\rm hom}\,(T^{n-1},F)\hookrightarrow {\rm
  hom}\,(T^n,F\otimes\Lambda^1)$ naturel en
l'objet $F$ de $\F_{n-1}$. C'est un
isomorphisme car l'application linéaire sous-jacente est inverse des
isomorphismes d'espaces vectoriels $$s_n(F\otimes\Lambda^1)={\rm
  hom}_\F (T^n,F\otimes\Lambda^1)\simeq {\rm hom}\,((T^n :
\Lambda^1),F)$$
$$\simeq {\rm hom}\,(T^{n-1},F)^{\oplus n}\simeq {\rm
  Ind}(s_{n-1}(F)).$$

Ainsi, le diagramme
\begin{equation*}
%\label{repdiv2b}
\xymatrix{\mathbf{Mod}_{\mathbb{F}_2[\Sigma_n]} & \mathbf{Mod}_{\mathbb{F}_2[\Sigma_{n-1}]}\ar[l]_-{{\rm Ind}} \\
\mathcal{F}_n\ar[u]^{s_n} & \mathcal{F}_{n-1}\ar[l]^-{-\otimes\Lambda^1}\ar[u]_{s_{n-1}}
}\end{equation*}
commute à isomorphisme naturel près.

La commutativité du diagramme (\ref{repdiv2}) s'en déduit  par adjonction.
\end{proof}

\begin{ex} En prenant $M=\mathbb{F}_2$, on obtient
  $(S^n:\Lambda^1)\simeq S^{n-1}$. Dualement, on a
  $\mathbf{Hom}\,(\Lambda^1,\Gamma^n)\simeq\Gamma^{n-1}$.

On peut montrer que $(\Gamma^n
  : \Lambda^1)\simeq\underset{i\in\mathbb{N}}{\bigoplus}
  \Gamma^{n-2^i}$, et donc
  $\mathbf{Hom}\,(\Lambda^1,S^n)\simeq\underset{i\in\mathbb{N}}{\bigoplus} S^{n-2^i}$ dualement ; cela illustre la nécessité de l'hypothèse d'homogénéité et de cohomogénité dans la proposition \ref{fonddivl} (\ref{pfo4}).
\end{ex}

\section{Compléments sur les foncteurs de Weyl}\label{cplw}

%\subsection{Quelques lemmes techniques}

Les considérations générales de la section précédente méritent d'être
appliquées aux foncteurs \go concrets\gf introduits au paragraphe
\ref{sws}, afin d'appréhender le comportement du foncteur de division par
$\Lambda^1$ en termes de facteurs de composition. C'est ce à quoi
s'emploie le paragraphe \ref{sdlw}.

Après des calculs préliminaires, nous donnons, au §\,\ref{spws}, un
résultat sur les facteurs de composition de certains foncteurs de Weyl
qui permettront, dans la section \ref{sctdl}, de leur appliquer efficacement des raisonnements
utilisant le foncteur $(-:\Lambda^1)$. Il s'agit d'éviter que le
facteur de composition $S_{\lambda^-_r}$ de $(W_\lambda:\Lambda^1)$
puisse déjà être facteur de composition de $({\rm rad}\,W_\lambda:\Lambda^1)$.

\subsection{Quelques lemmes techniques}

Nous commençons par donner la définition et les propriétés de base de morphismes qui joueront un rôle
important dans le paragraphe suivant.

\begin{nota}\label{proj-palt}
Pour $i\geq j\geq 1$, nous noterons $\Pi_{i,j}$ et $\Pi'_{i,j}$ les
endomorphismes de $\Lambda^i\otimes\Lambda^j$ donnés respectivement
par les compositions suivantes.
$$\Pi_{i,j}=\Lambda^i\otimes\Lambda^j\xrightarrow{D\theta}\Lambda^{i-1}\otimes\Lambda^{j+1}\xrightarrow{\theta}\Lambda^i\otimes\Lambda^j$$
$$\Pi'_{i,j}=\Lambda^i\otimes\Lambda^j\xrightarrow{\theta}\Lambda^{i+1}\otimes\Lambda^{j-1}\xrightarrow{D\theta}\Lambda^i\otimes\Lambda^j$$
\end{nota}

\begin{lm}
\label{calt}
\begin{enumerate}
\item \label{p2} Pour $k+l\leq j$, on a $$\theta_{i+k,j-k,l}\theta_{i,j,k}=\frac{(k+l)!}{k!\,l!}\cdot\theta_{i,j,k+l}\,.$$
\item \label{lbproj} Pour $i\geq j\geq 1$, on a $\quad\Pi_{i,j}+\Pi'_{i,j}=(i+j)\,id$.
\item \label{comt} Soient $i, j, k, t, u$ des entiers positifs, avec $j\geq t+u$. Le diagramme suivant est commutatif.
$$\xymatrix{\Lambda^{i,j,k} \ar[r]^-{\theta\otimes id}\ar[d]_{id\otimes D\theta} & \Lambda^{i+t,j-t,k}\ar[d]^{id\otimes D\theta} \\
\Lambda^{i,j-u,k+u}\ar[r]^-{\theta\otimes id} & \Lambda^{i+t,j-t-u,k+u}}$$
\end{enumerate}
\end{lm}

\begin{proof} Ces calculs sont analogues ; établissons par exemple
  (\ref{lbproj}). Soient $V$ un espace vectoriel de dimension
  finie et
  $a_1,\dots,a_i ; b_1,\dots, b_j$ des éléments de $V$. Si
  $P=\{t_1<\dots<t_k\}$ est une partie de $\mathbf{i}=\{1,\dots,i\}$,
  notons $a^{\wedge P}$ pour $a_{t_1}\wedge\dots\wedge a_{t_k}$. On a 
$$D\theta_{i-1,j+1,1}(V)(a^{\wedge\mathbf{i}}\otimes b^{\wedge\mathbf{j}})=\sum_{k=1}^i a^{\wedge \{k\}^c}\otimes (b^{\wedge\mathbf{j}}\wedge a_k)$$
 où l'exposant $c$ indique le complémentaire ensembliste ; puis
$$\Pi_{i,j}(V)(a^{\wedge\mathbf{i}}\otimes b^{\wedge\mathbf{j}})=\sum_{k=1}^i \Big(a^{\wedge\mathbf{i}}\otimes b^{\wedge\mathbf{j}}+\sum_{l=1}^j (a^{\wedge \{k\}^c}\wedge b_l)\otimes (b^{\wedge\{l\}^c}\wedge a_k)\Big).$$

De même
$$\Pi'_{i,j}(V)(a^{\wedge\mathbf{i}}\otimes b^{\wedge\mathbf{j}})=\sum_{l=1}^j \Big(a^{\wedge\mathbf{i}}\otimes b^{\wedge\mathbf{j}}+\sum_{k=1}^i (a^{\wedge \{k\}^c}\wedge b_l)\otimes (b^{\wedge\{l\}^c}\wedge a_k)\Big),$$
d'où l'assertion (\ref{lbproj}).\end{proof}

\begin{pr}
\label{pproj}
Soient $i\geq j\geq 1$ des entiers.
\begin{enumerate}
\item \label{proj1} Supposons $i-j$ impair : 
\begin{enumerate}
\item \label{proj11} $\Pi_{i,j}$ et $\Pi'_{i,j}$ sont deux projecteurs
  dont la somme est l'identité ;
\item \label{proj12}
  $im\,\Pi_{i,j}=ker\,\Pi'_{i,j}=im\,\theta_{i-1,j+1,1}=ker\,\theta_{i,j,1}\supset W_{(i,j)}$ ;
\item \label{proj13} $im\,\Pi'_{i,j}=ker\,\Pi_{i,j}=im\,D\theta_{i,j,1}=ker\,D\theta_{i-1,j+1,1}$.
\end{enumerate} 
\item \label{proj2} Supposons $i-j$ pair. Alors $\Pi_{i,j}=\Pi'_{i,j}$ a une image dont les facteurs de composition sont du type $S_{(i+t,j-t)}$ avec $1\leq t\leq j$.
\end{enumerate}
\end{pr}

\begin{proof} Le point (\ref{p2}) du lemme \ref{calt} montre que
  l'on a toujours $\Pi'_{i,j}\Pi_{i,j}=0$. Combiné avec le point (\ref{lbproj}) de ce même lemme, ce fait implique (\ref{proj11}).

On a ensuite $im\,\Pi_{i,j}\subset im\,\theta_{i-1,j+1,1}$, et
$ker\,\Pi'_{i,j}\supset ker\,\theta_{i,j,1}\supset W_{(i,j)}$, d'où
l'on déduit (\ref{proj12}) via (\ref{proj11}) et l'égalité
$ker\,\theta_{i,j,1}=im\,\theta_{i-1,j+1,1}$ (valable car $i>j$)
déduite de la proposition 1.3.1 de \cite{Franjou}.

L'assertion (\ref{proj13}) se prouve de façon analogue.

\smallskip

Pour le point (\ref{proj2}), on utilise le point (\ref{lbproj}) du lemme \ref{calt} et le fait que les facteurs de composition de $\Lambda^{i+1}\otimes\Lambda^{j-1}$, dont $im\,\Pi'_{i,j}$ est un sous-quotient, sont les $S_{(i+t,j-t)}$ pour $1\leq t\leq j$.\end{proof}

Les autres résultats techniques dont nous aurons besoin se démontrent
très facilement à l'aide de la notion classique suivante. Ils
interviendront dans la section \ref{sec4}.

\begin{defi}
Soient $V\in {\rm Ob}\,\E^f$ et $\lambda$ une partition de longueur
$r$.
\begin{enumerate}\item Soient $v_1,\dots,v_{\lambda_1}$ des éléménts
  de $V$. L'{\em élément
  semi-standard associé à $\lambda$ et
  $v_1,\dots,v_{\lambda_1}$} est l'élément de $\Lambda^\lambda(V)$ défini par
$$s^{st}_\lambda(v_1,\dots,v_{\lambda_1})=\bigotimes_{i=1}^r
(v_1\wedge\dots\wedge v_{\lambda_i}).$$
\item Soient $a_{i,j}$ ($1\leq i\leq r$, $1\leq j\leq\lambda_i$) des éléments de $V$. L'{\em élément
  standard associé à $\lambda$ et
  $(a_{i,j})$}  est l'élément de $\Lambda^\lambda(V)$ défini par
$$g^{st}\big((a_{i,j})\big)=\sum_{\sigma\in R_{\lambda}}\bigotimes_{i=1}^r
\big(a_{\sigma(i,1)}\wedge\dots\wedge a_{\sigma(i,\lambda_i)}\big)$$ 
où $R_{\lambda}$ désigne le groupe des permutations de l'ensemble
$\{(i,j)\,|\,1\leq i\leq r,\:1\leq j\leq\lambda_i\}$ laissant
invariante la deuxième composante. 
\end{enumerate}
\end{defi}

\begin{pr}\label{semist}
Soient $V\in {\rm Ob}\,\E^f$ et $\lambda$ une partition.
\begin{enumerate}\item L'espace vectoriel $W_\lambda(V)$ est le sous-espace de $\Lambda^\lambda(V)$ engendré par les éléments semi-standard $s^{st}_\lambda(v_1,\dots,v_{\lambda_1})$ pour $v_1,\dots,v_{\lambda_1}\in V$.
\item Si $\lambda$ est régulière, $W_\lambda(V)$ est le sous-espace de $\Lambda^\lambda(V)$ engendré par les éléments standard.
\end{enumerate}
\end{pr}

Cette propriété, pour laquelle  nous renvoyons à \cite{James}
(cf. aussi \cite{PS}, §\,$2$, pour le cas de $\F$), permet de
simplifier de façon appréciable certains calculs sur les foncteurs de
Weyl.

\begin{lm}
\label{mweyl}
Soient $i$, $j$, $t$ des entiers naturels tels que $i\geq j$ ; notons $\tau : \Lambda^i\otimes\Lambda^j\to\Lambda^j\otimes\Lambda^i$ l'isomorphisme d'échange des deux facteurs du produit tensoriel. La restriction à $W_{(i,j)}$ de $\Lambda^i\otimes\Lambda^j\xrightarrow{\tau}\Lambda^j\otimes\Lambda^i\xrightarrow{\theta_{j,i,i-j}}\Lambda^i\otimes\Lambda^j$ coïncide avec l'identité.
\end{lm}

\begin{proof} Un calcul immédiat montre en effet que le morphisme
  $\theta\circ\tau$ ne modifie pas les éléments semi-standard. 
\end{proof}

\begin{pr}\label{prpra}
Soient $r\in\mathbb{N}^*$ et $\lambda$ une partition régulière de
longueur $r$ ou $r-1$. Notons $j_{\lambda,r}$ le morphisme
$$\Lambda^\lambda\otimes\Lambda^r\hookrightarrow\Lambda^\lambda\otimes T^r\simeq\underset{i=1}{\bigotimes^r}(\Lambda^{\lambda_i}\otimes\Lambda^1)\xrightarrow{\underset{i}{\bigotimes}\theta_{\lambda_i,1,1}}\Lambda^{(\lambda_1+1,\dots,\lambda_r+1)}\,.$$

On a $j_{\lambda,r}(W_\lambda\otimes\Lambda^r)=W_{(\lambda_1+1,\dots,\lambda_r+1)}$.
\end{pr}

\begin{proof} Cela résulte de la proposition \ref{semist} et du calcul suivant : si les $a_{i,j}$ ($1\leq i\leq r$, $1\leq j\leq\lambda_i+1$) sont des éléments d'un espace vectoriel $V$, on a 
$j_{\lambda,r}(V)\big(g^{st}\big((a_{i,j})_{j\leq\lambda_i}\big)\otimes (a_{1,\lambda_1+1}\wedge\dots\wedge a_{r,\lambda_r+1})\big)=g^{st}\big((a_{i,j})\big)\,.$
\end{proof}

\begin{nota}
\label{crweyl2}
Soient $\lambda_1>\dots>\lambda_r\geq 0$ et $n\geq 0$ des entiers. On
pose $\lambda_{+i}=(\lambda_1+i\dots,\lambda_r+i)$. On définit par
récurrence une flèche\\
$\Lambda^\lambda\otimes T^n(\Lambda^r)\xrightarrow{j^n_\lambda}\Lambda^{\lambda_{+n}}$ par $j^0_\lambda=id$ et
$$j^n_\lambda=\Lambda^\lambda\otimes\Lambda^r\otimes T^{n-1}(\Lambda^r)\xrightarrow{j_{\lambda,r}\otimes T^{n-1}(\Lambda^r)}\Lambda^{\lambda_{+1}}\otimes  T^{n-1}(\Lambda^r)\xrightarrow{j^{n-1}_{\lambda_{+1}}}\Lambda^{\lambda_{+n}}$$
si $n>0$. 
\end{nota}

On déduit de la proposition \ref{prpra} le résultat suivant.

\begin{cor}
\label{crweyl} Soient $\lambda_1>\dots>\lambda_r\geq 0$ et $n\geq 0$ des entiers.
\begin{itemize}
\item Il existe un unique morphisme $w^{\lambda_1,\dots,\lambda_r}_n : \Lambda^\lambda\otimes\Lambda^n(\Lambda^r)\to\Lambda^{\lambda_{+n}}$ tel que le diagramme suivant commute.
$$\xymatrix{\Lambda^\lambda\otimes T^n(\Lambda^r) \ar[r]^(.62){j^n_\lambda}\ar@{>>}[d] & \Lambda^{\lambda_{+n}} \\
\Lambda^\lambda\otimes\Lambda^n(\Lambda^r)\ar@{-->}[ur]_(.57){w^{\lambda_1,\dots,\lambda_r}_n} &
}$$
\item On a $w^{\lambda_1,\dots,\lambda_r}_n\big(W_\lambda\otimes\Lambda^n(\Lambda^r)\big)=W_{\lambda_{+n}}$.
\item Pour $i>j>t>0$, le diagramme suivant commute.
$$\xymatrix{\Lambda^i\otimes\Lambda^j\otimes\Lambda^n(\Lambda^3) \ar[r]^-{\theta\otimes id} \ar[d]_{w^{i,j,0}_n} & \Lambda^{i+t}\otimes\Lambda^{j-t}\otimes\Lambda^n(\Lambda^3)\ar[d]^{w^{i+t,j-t,0}_n} \\
\Lambda^{i+n}\otimes\Lambda^{j+n}\otimes\Lambda^n\ar[r]^-{\theta\otimes id} & \Lambda^{i+t+n}\otimes\Lambda^{j-t+n}\otimes\Lambda^n 
}$$
\end{itemize}
\end{cor}

\subsection{Partitions  Weyl-séparantes}\label{spws}

\begin{nota}\label{notpl2}
Soit $\lambda=(\lambda_1,\dots,\lambda_r)$ un $r$-uplet d'entiers
et $1\leq i\leq r$. On définit $\lambda^-_i$ (resp. $\lambda^+_i$) par
$(\lambda^-_i)_j=\lambda_j$ si $j\neq i$ et
$(\lambda^-_i)_i=\lambda_i-1$ (resp. $(\lambda^+_i)_j=\lambda_j$ si
$j\neq i$ et $(\lambda^+_i)_i=\lambda_i+1$).

Pour alléger ces notations, nous simplifierons des écritures du type $(\lambda^+_a)^-_b$ en $\lambda^{+,-}_{a,b}$.
\end{nota}

\begin{defi} On dit
  qu'une partition $\lambda$ de longueur $r$ est :
\begin{itemize}
\item
  \textbf{Weyl-séparante} (ou \textbf{W-séparante}) s'il
  n'existe pas de partition régulière $\mu$ de $|\lambda|$ telle que
  $\mu\leq\lambda^{+,-}_{1,r}$ et $\mu\vdash {\rm rad}\,W_\lambda$,
\item \textbf{alternée} si $\lambda_i-\lambda_{i+1}$ est impair pour $1\leq i<r$.
\end{itemize}
\end{defi}

On notera qu'une partition alternée est toujours régulière.

\smallskip

Intuitivement, une partition régulière $\lambda$ est Weyl-séparante si le
foncteur de Weyl associé n'a pas de facteur de composition \go proche\gf mais
distinct de $S_\lambda$. Cette notion sera fort utile pour détecter
certains facteurs de composition, comme nous le verrons dans la
section \ref{sctdl}.

Le but de ce paragraphe est d'établir qu'une partition alternée est
Weyl-séparante. Cette propriété, analogue aux considérations de
\cite{James}, §\,$24$, fournira tous les cas  de
Weyl-séparation dont nous aurons besoin.

\begin{lm}
\label{wsf}
Soit $\lambda$ une partition régulière de longueur $r$. La
partition $\lambda$ est W-séparante si et seulement s'il n'existe pas de
partition régulière $\mu$ de $|\lambda|$ telle que
$\lambda\leq\mu\leq\lambda^{+,-}_{1,r}$ et
$$\mu\vdash W_\lambda\cap\underset{1\leq i\leq
  r-1}{\sum}im\,D\psi^{i,1}_\lambda$$ où l'on emploie le morphisme $\psi^{i,t}_\lambda$ de la notation \ref{nfd}.
\end{lm}

\begin{proof} On utilise le théorème \ref{repf} (\ref{nrf}), en notant que si $\mu$ est une
partition régulière de $|\lambda|$ telle que $\mu\vdash
\Lambda^{\lambda_1,\dots,\lambda_{i-1},\lambda_i+t,\lambda_{i+1}-t,\lambda_{i+2},\dots,\lambda_r}$  
pour un  $1\leq i\leq r-1$ et un $2\leq t\leq\lambda_{i+1}$, alors $\mu$ n'est pas inférieure à $\lambda^{+,-}_{1,r}$ grâce au théorème \ref{fclambda}.\end{proof}

\begin{nota} Soit $L$ une $\mathbb{F}_2$-algèbre associative et
  unitaire. Nous noterons, dans ce paragraphe, $[a,b]=aba+bab$ pour $(a,b)\in L^2$. 
\end{nota}

\begin{rem} Il s'agit d'une notation {\em ad hoc} valable uniquement
  dans ce paragraphe. Elle est motivée par le fait que ce \go crochet\gf joue intuitivement le même rôle qu'un commutateur usuel sur les endomorphismes que l'on considère dans la démonstration de la proposition \ref{pfc} ci-après.
\end{rem}

\begin{lm}\label{lmbiz}
\begin{enumerate}
\item \label{biz1} Si $a$ et $b$ sont deux idempotents de $L$, on a $[1+a,1+b]=[a,b]$.
\item \label{biz2} Soient $k\in\mathbb{N}^*$, $u_0,\dots,u_k$ des
  éléments de $L$ tels que $u_i u_j=u_j u_i$ si $|i-j|\geq 2$,
  $v=u_k\dots u_1$ et $I$ l'idéal bilatère de $L$ engendré par les
  $[u_{i-1},u_i]$ ($1\leq i\leq k$). Alors $vu_0 v\in I+Lu_0$.
\end{enumerate}
\end{lm}

\begin{proof}  L'assertion (\ref{biz1}) provient du calcul suivant
  :
$$(1+a)(1+b)(1+a)=1+a+b+ab+ba+aba,$$ et $(1+b)(1+a)(1+b)=1+a+b+ab+ba+bab$.

Pour l'assertion (\ref{biz2}), on raisonne par récurrence sur
$k$. Le cas $k=1$ provient de ce que $u_1 u_0 u_1=[u_0,u_1]+u_0 u_1
u_0$, le cas $k=2$ de ce que $$u_2 u_1 u_0 u_2 u_1=u_2 u_1 u_2 u_0
u_1=[u_1,u_2] u_0 u_1 + u_1 u_2 u_1 u_0 u_1$$
$$=[u_1,u_2] u_0 u_1 + u_1 u_2 [u_0,u_1] + u_1 u_2 u_0 u_1 u_0.$$ 

On a donc, en supposant maintenant $k>2$,
$$u_k\dots u_1 u_0 u_k\dots u_1=(u_k\dots u_2)u_1 u_0 (u_k\dots u_2)
u_1$$
$$=[u_1,u_k\dots u_2] u_0 u_1 + u_1 u_k\dots u_2 [u_0,u_1] + u_1 u_k\dots u_2 u_0 u_1 u_0\,.$$

Or l'hypothèse de récurrence montre que $[u_1,u_k\dots u_2]\in I+ L
u_1$. Par conséquent,
$$u_k\dots u_1 u_0 u_k\dots u_1\in I+ L u_1 u_0 u_1 + L u_0\,,$$
ce qui permet de conclure grâce au cas $k=1$ (qui montre que $u_1 u_0 u_1\in I+L u_0$).\end{proof}

\begin{lm}\label{lfdc}
Etant donnée une partition alternée $(i,j,k)$ de longueur $3$, notons $A$ (resp. $B$) le projecteur (cf. proposition \ref{pproj}) $\Pi_{i,j}\otimes\Lambda^k$ (resp. $\Lambda^i\otimes\Pi_{j,k}$) de $\Lambda^i\otimes\Lambda^j\otimes\Lambda^k$. Avec la notation du lemme précédent, l'image de $[A,B]$ n'a pas de facteur de composition $S_\lambda$ si $\lambda$ est une partition régulière de $i+j+k$ telle que $\lambda\leq (i+1,j,k-1)$.
\end{lm}

\begin{proof} D'après l'assertion (\ref{biz1}) du lemme précédent,
  $[A,B]=[A',B']$ où $A'=\Pi'_{i,j}\otimes\Lambda^k$ et
  $B'=\Lambda^i\otimes\Pi'_{j,k}$. Il suffit donc de montrer que les
  deux endomorphismes $A'B'A'$ et $B'A'B'$ ont une image sans facteur
  de composition du type indiqué dans l'énoncé. Considérons pour cela
  le diagramme suivant.
$$\xymatrix{\Lambda^{i,j,k}\ar[r]^-{\theta\otimes id} \ar[dr]_{A'} & \Lambda^{i+1,j-1,k} \ar[d]^-{D\theta\otimes id}\ar[r]^-{id\otimes\theta} & \Lambda^{i+1,j,k-1} \ar[d]_(.47){D\theta\otimes id} \ar[dr]^{\Pi_{i+1,j}\otimes id} & \\
& \Lambda^{i,j,k} \ar[r]^(.46){id\otimes\theta}\ar[dr]_{B'} & \Lambda^{i,j+1,k-1} \ar[r]_(.47){\theta\otimes id}\ar[d]^{id\otimes D\theta} & \Lambda^{i+1,j,k-1}\ar[d]^{id\otimes D\theta} \\
& & \Lambda^{i,j,k}\ar[r]^(.46){\theta\otimes id}\ar[dr]_{A'} & \Lambda^{i+1,j-1,k}\ar[d]^{D\theta\otimes id} \\
& & & \Lambda^{i,j,k}
}$$
 
Il commute (les triangles, par définition des flèches obliques, et les carrés par la propriété (\ref{comt}) du lemme \ref{calt}), donc l'image de $A'B'A'$ est un sous-quotient de l'image de $\Pi_{i+1,j}\otimes\Lambda^{k-1}$.

On établit de même que l'image de $B'A'B'$ est un sous-quotient de l'image de $\Lambda^{i+1}\otimes\Pi_{j,k-1}$.

Pour conclure, on utilise le point (\ref{proj2}) de la
proposition \ref{pproj} et le théorème \ref{fclambda} : ils montrent que les facteurs de composition de degré $i+j+k$ de ces images sont associés à des partitions supérieures à $(i+2,j-1,k-1)$ dans le premier cas, et à $(i+1,j+1,k-2)$ dans le second.\end{proof}

\begin{pr}
\label{pfc}
Une partition alternée est W-séparante.
\end{pr}

\begin{proof} Soit $\lambda$ une partition alternée de longueur $r$. On définit des endomorphismes de $\Lambda^\lambda$ par
$$p_i=\Lambda^{\lambda_1,\dots,\lambda_{i-1}}\otimes\Pi_{\lambda_i,\lambda_{i+1}}\otimes\Lambda^{\lambda_{i+2},\dots,\lambda_r}\;\,(1\leq i\leq r-1)\;\,\text{et}\; P=\prod_{i=1}^{r-1}\left(\prod_{j=1}^{r-i}p_{r-j}\right).$$

Par l'assertion (\ref{proj13}) de la proposition \ref{pproj},
$im\,D\psi_{i,1}=im\,(1+p_i)$. Comme l'image de chacun des projecteurs $p_i$ contient $W_\lambda$ (utiliser l'assertion (\ref{proj12}) de la proposition \ref{pproj}), on en déduit
$$W_\lambda\cap\underset{1\leq i\leq r-1}{\sum}im\,D\psi_{i,1}\subset\underset{1\leq i\leq r-1}{\sum}im\,P(1+p_i)\,.$$

Il suffit donc d'établir, grâce au lemme \ref{wsf}, que $S_\mu$ n'est
facteur de composition d'aucune des images des $P(1+p_i)$ si  $\mu$
est une
partition régulière de $|\lambda|$ telle que $\lambda\leq\mu\leq\lambda^{+,-}_{1,r}$. Pour cela, on note que $p_i$ et $p_j$ commutent si $|i-j|\geq 2$, de sorte que $P(1+p_i)$ appartient à l'idéal bilatère de ${\rm End}\,\Lambda^\lambda$ engendré par $p_{r-1}\dots p_{i+1}p_i p_{r-1}\dots p_{i+1}(1+p_i)$. On applique ensuite l'assertion (\ref{biz2}) du lemme \ref{lmbiz} pour obtenir que $P(1+p_i)$ appartient à l'idéal bilatère engendré par les $[p_j,p_{j+1}]$, puisque $p_i(1+p_i)=0$. Le lemme \ref{lfdc}, combiné au théorème \ref{fclambda}, permet alors de conclure.\end{proof}

\subsection{Division par $\Lambda^1$ des foncteurs de Weyl}\label{sdlw}
%\subsection{Calcul sur les foncteurs de Weyl}

Nous donnons dans cette sous-section quelques résultats de base sur la division par $\Lambda^1$ des foncteurs de Weyl.  Ces résultats
joueront un rôle fondamental pour la détection de facteurs de
composition à l'aide du foncteur $(\cdot : \Lambda^1)$ ; de fait, on ne dispose pratiquement d'aucun autre renseignement sur la division par $\Lambda^1$ d'un foncteur simple que ceux que l'on déduit grossièrement du cas des foncteurs de Weyl.

\begin{pr}\label{fcsdl}Soit $\lambda$ une partition régulière de
  longueur $r$ de $n$. 
\begin{enumerate}\item Les morphismes $\upsilon_{W_\lambda} :
  \mathbf{Hom}(\Lambda^1,W_\lambda)\to (W_\lambda : \Lambda^1)$ et \\
$\upsilon_{S_\lambda} :
  \mathbf{Hom}(\Lambda^1,S_\lambda)\to (S_\lambda : \Lambda^1)$
  sont des isomorphismes.
\item Si $\mu$ est une partition régulière telle que $\mu\vdash (W_{\lambda}:\Lambda^1)$, on a :
\begin{itemize}
\item soit $|\mu|=n-1$ et $\mu\geq\lambda^-_{1}$,
\item soit $|\mu|<n-1$, $l(\mu)<r$, $\mu_{1}\geq\lambda_{1}$ et $\mu_{r-1}\leq\lambda_{r}$. 
\end{itemize}
\end{enumerate}
\end{pr}

\begin{proof} Le foncteur $W_{\lambda}$ est homogène car inclus dans
  $\Lambda^{\lambda}$, et cohomogène par le théorème/définition \ref{repf}. La
  proposition \ref{fonddivl} fournit donc le premier point. On en
  déduit
$$(S_\lambda : \Lambda^1)\twoheadleftarrow (W_\lambda :
  \Lambda^1)\simeq\mathbf{Hom}(\Lambda^1,W_\lambda)\hookrightarrow\mathbf{Hom}(\Lambda^1,\Lambda^\lambda).$$
On conclut en appliquant le corollaire \ref{cordiv} et le théorème \ref{fclambda}.\end{proof}

\begin{rem} Le morphisme $\upsilon_{W_\lambda}$ n'est pas forcément un
  isomorphisme lorsque $\lambda$ est une partition non régulière
  (cf. remarque \ref{rqnsw}).
\end{rem}

\begin{rem} L'hypothèse de la seconde assertion de la proposition
  \ref{fcsdl} est en particulier satisfaite si $\mu\vdash (S_{\lambda}:\Lambda^1)$.
\end{rem}

Précisons cette proposition par l'analogue suivant du théorème de branchement de James pour la restriction des modules de Specht en théorie des représentations du groupe symétrique (cf. \cite{James}) ; le seul ingrédient nouveau par rapport à la théorie des représentations est la proposition \ref{fonddivl}. Nous donnons une démonstration directe fondée sur le lemme calculatoire simple suivant.

\begin{lm}\label{lmcd}
Soient $i$, $j$, $t$ des entiers positifs tels que $t\leq j$. Le morphisme
$$(\Lambda^{i-1}\otimes\Lambda^j)\oplus(\Lambda^i\otimes\Lambda^{j-1})\simeq
(\Lambda^i\otimes\Lambda^j:\Lambda^1)\xrightarrow{(\theta_{i,j,t}:\Lambda^1)}(\Lambda^{i+t}\otimes\Lambda^{j-t}:\Lambda^1)$$
$$\simeq(\Lambda^{i+t-1}\otimes\Lambda^{j-t})\oplus(\Lambda^{i+t}\otimes\Lambda^{j-t-1})$$
vérifie les propriétés suivantes :
\begin{itemize}
\item sa composante $\Lambda^{i-1}\otimes\Lambda^j\to\Lambda^{i+t-1}\otimes\Lambda^{j-t}$ est égale à $\theta_{i-1,j,t}$,
\item sa composante $\Lambda^{i-1}\otimes\Lambda^j\to\Lambda^{i+t}\otimes\Lambda^{j-t-1}$ est nulle,
\item sa composante $\Lambda^i\otimes\Lambda^{j-1}\to\Lambda^{i+t}\otimes\Lambda^{j-t-1}$ est égale à $\theta_{i,j-1,t}$.
\end{itemize}  
\end{lm}

\begin{proof} Par adjonction, il s'agit de vérifier la commutativité des diagrammes suivants :
$$\xymatrix{\Lambda^1\otimes\Lambda^{i-1}\otimes\Lambda^j\ar[r]^-{D\theta\otimes id}\ar[d]_{\Lambda^1\otimes\theta} & \Lambda^i\otimes\Lambda^j\ar[d]^\theta \\
\Lambda^1\otimes\Lambda^{i+t-1}\otimes\Lambda^{j-t}\ar[r]^-{D\theta\otimes id} & \Lambda^{i+t}\otimes\Lambda^{j-t}
}$$
pour les deux premières composantes considérées, car on peut remplacer
$(-:\Lambda^1)$ par $\mathbf{Hom}(\Lambda^1,-)$ grâce à la proposition \ref{fcsdl}, et
$$\xymatrix{\Lambda^i\otimes\Lambda^j\ar[d]_\theta\ar[r]^-{id\otimes D\theta} & \Lambda^i\otimes\Lambda^{j-1}\otimes\Lambda^1\ar[d]^{\theta\otimes\Lambda^1}\\
\Lambda^{i+t}\otimes\Lambda^{j-t}\ar[r]^(.4){id\otimes D\theta} & \Lambda^{i+t}\otimes\Lambda^{j-t-1}\otimes\Lambda^1
}$$
pour la dernière.

On conclut maintenant grâce à l'assertion \ref{comt} du lemme \ref{calt}.\end{proof}

\begin{pr}
\label{cbdiv} Soit $\lambda$ une partition régulière de longueur $r$. Il existe une filtration
$0=F_{0}\subset F_{1}\subset\dots\subset
F_{r}=(W_{\lambda}:\Lambda^1)$ telle que
$F_{i}/F_{i-1}\simeq W_{\lambda^-_{i}}$ pour $0<i\leq r\,$.

Précisément, $F_i$ est donné par le diagramme commutatif cartésien d'inclusions
\begin{equation*}
\xymatrix{F_i\ar@{^{(}->}[r] \ar@{^{(}->}[d] \ar@{}[dr] | \square & \underset{1\leq
    j\leq i}{\bigoplus} \Lambda^{\lambda^-_j}\ar@{^{(}->}[d] \\
(W_\lambda:\Lambda^1) \ar@{^{(}->}[r] & (\Lambda^\lambda :\Lambda^1).
}\end{equation*}
%
%En particulier, la composée $(W_\lambda:\Lambda^1)\hookrightarrow (\Lambda^\lambda:\Lambda^1)\twoheadrightarrow\Lambda^{\lambda^-_r}$ a pour image $W_{\lambda^-_r}$. L'épimorphisme induit 
%$\pi^W_\lambda : (W_\lambda:\Lambda^1)\twoheadrightarrow W_{\lambda^-_r}$  peut se voir comme l'adjoint de l'inclusion $\gamma^W_\lambda : W_\lambda\hookrightarrow W_{\lambda^-_r}\otimes\Lambda^1$.
\end{pr}

\begin{proof} La proposition \ref{fcsdl} permettant de remplacer $(-:\Lambda^1)$ par $\mathbf{Hom}(\Lambda^1,-)$, qui est exact à gauche, le lemme précédent permet de conclure en utilisant la remarque \ref{remweyl}.\end{proof}

\begin{rem}\label{rqnsw}
La proposition est en défaut pour une partition non régulière, en
raison de la non cohomogénéité du foncteur de Weyl associé. Par exemple, $(W_{1,1}:\Lambda^1)=(\Gamma^2:\Lambda^1)\simeq\mathbb{F}_2\oplus\Lambda^1$.
\end{rem}

\section{Détection de facteurs de composition par division par $\Lambda^1$}\label{sctdl}

%\subsection{Préliminaires formels}\label{sectf}

  Le problème de l'effet de foncteurs remarquables sur les facteurs de
composition d'objets d'une catégorie abélienne se rencontre
naturellement dans divers contextes, l'étude directe des facteurs de
composition s'avérant généralement ardue, voire inabordable. Le cas le
plus simple d'un foncteur exact, et commutant aux colimites si l'on
s'intéresse à des objets seulement localement finis, se révèle souvent
insuffisant.

Dans \cite{Po2}, Powell introduit des endofoncteurs $\tilde{\nabla}_n$ de la catégorie
$\F$ qui ne sont exacts ni
à gauche ni à droite, mais qui préservent les épimorphismes et les
monomorphismes. Ainsi, si le simple $S$ est facteur de composition de
$F$, alors $\tilde{\nabla}_n S$ est un sous-quotient de
$\tilde{\nabla}_n F$. Comme pour certains foncteurs simples $S$,
$\tilde{\nabla}_n S$ est un foncteur simple explicite, on obtient des
renseignements sur les facteurs de composition de $\tilde{\nabla}_n F$
dès lors que l'on connaît certains facteurs de composition de
$F$. Powell a montré le grand intérêt des foncteurs
$\tilde{\nabla}_n$ dans ses articles \cite{Po1} et \cite{Po3}.

Notre démarche, utilisant le foncteur $(-:\Lambda^1)$,
est inverse : nous cherchons à obtenir des renseignements sur les
facteurs de composition d'un foncteur analytique $F$, connaissant certains facteurs de
composition de $(F:\Lambda^1)$. Ceci est théoriquement possible dans
la mesure où si $\lambda\vdash (F:\Lambda^1)$, alors il existe $\mu$
telle que $\mu\vdash F$ et $\lambda\vdash (S_\mu:\Lambda^1)$, puisque
le foncteur $(-:\Lambda^1)$ commute aux colimites. Il semble cependant
illusoire d'obtenir des résultats très généraux, en raison des deux
écueils suivants :
\begin{enumerate}\item la description des facteurs de composition de
  $(S_\mu:\Lambda^1)$ est hors de portée en général ;
\item un foncteur simple est en général facteur de composition de la
  division par $\Lambda^1$ d'un grand nombre de foncteurs simples.
\end{enumerate}

Afin de contourner ces difficultés, nous mettons deux restrictions à
notre problème initial :
\begin{enumerate}\item on suppose que le foncteur $F$ est un
  sous-objet d'un foncteur connu $X$, et l'on cherche s'il contient ou non
 des facteurs de composition identifiés dans $X$ ; 
\item on se limite au cas où l'on maîtrise un tant soit peu l'effet de la division par
  $\Lambda^1$ sur le foncteur simple que l'on cherche à détecter.
\end{enumerate}

La stratégie de détection dans un sous-objet d'un objet connu est
présentée dans un cadre général dans le paragraphe \ref{sectf}. Dans
la catégorie $\F$, la
seconde restriction nous amènera à travailler sur les partitions
Weyl-séparantes, qui ont été introduites à cette fin, comme nous le
verrons au paragraphe \ref{dfgm}. Nous terminons cette section avec un
résultat technique plus global, la proposition \ref{cr-fd1}, adapté à la situation des foncteurs
$\bar{I}^{\otimes 2}\otimes\Lambda^n$ que nous avons en vue. 

\subsection{Préliminaires formels}\label{sectf}

\begin{conv}
Dans ce paragraphe, $\mathcal{A}$ et $\mathcal{B}$ sont deux
catégories abéliennes avec colimites filtrantes exactes, $\Phi :
\mathcal{A}\to\mathcal{B}$ est un foncteur commutant aux colimites (en
particulier, exact à droite), $S$ (resp. $S'$) un objet simple de $\mathcal{A}$ (resp. $\mathcal{B}$).
\end{conv}

\begin{rem} L'intérêt d'avoir affaire à des catégories avec colimites filtrantes exactes pour
le maniement des facteurs de composition provient de l'observation
suivante : si un objet $X$ d'une telle catégorie est colimite
filtrante de sous-objets $A_i$, un objet simple $S$ est facteur de
composition de $X$ si et seulement s'il est facteur de composition de
l'un des $A_i$. On le voit en utilisant l'isomorphisme canonique
$Y\simeq\underset{i}{\col} (A_i\cap Y)$
valable pour tout sous-objet $Y$ de $X$ sous l'hypothèse d'exactitude
des colimites filtrantes (cf. \cite{Gab}).
\end{rem}

\begin{defi}\label{dedet}
On dit que l'objet simple $S$ est \textbf{$\Phi$-détecté par
  $S'$} dans un objet $X$
de $\mathcal{A}$  si $S'$ est facteur de
composition de $\Phi X$ et que  $S$ est facteur de composition de tout sous-objet $A$ de $X$ tel que $S'$ est
facteur de composition de $im\,\Phi(A\hookrightarrow X)$.

Plus généralement, si $\Phi X\xrightarrow{\pi} B$ est une flèche de
$\mathcal{B}$, on dit que $S$ est  \textbf{$\Phi$-détecté par $S'$
  dans $X$ relativement à $\pi$} si $S'$ est facteur de composition de
$im\,\pi$ et que $S$ est facteur de composition de tout sous-objet $A$ de $X$ tel que $S'$ est
facteur de composition de $\pi\big(im\,\Phi(A\hookrightarrow X)\big)$.
\end{defi}

\begin{rem}\label{remf}\begin{enumerate} 
\item Par exactitude à droite de $\Phi$, la suite
$$0\to im\,\Phi(A\hookrightarrow X)\to\Phi X\to\Phi(X/A)\to 0$$
est exacte. En particulier, l'hypothèse sur l'image est satisfaite si $S'$ n'est pas facteur de composition de $\Phi(X/A)$.
\item L'objet simple $S$ est $\Phi$-détecté par $S'$
  dans $S$ relativement à $\pi$ si et seulement si $S'$ est facteur de
  composition de $\pi(\Phi S)$.
\item La notion de $\Phi$-détection relativement à $\pi$ ne dépend que
  de $ker\,\pi$ : composer $\pi$ à gauche par un monomorphisme ne
  change pas la notion obtenue. 
\end{enumerate}
\end{rem}

\begin{lm}
\label{p-tech}
Soit $X$ un objet localement fini de $\mathcal{A}$. Les assertions
suivantes sont équivalentes.
\begin{enumerate}
\item Si un objet simple $T$ de $\mathcal{A}$ est facteur de composition de $X$, $S'$ n'est pas facteur de composition de $\Phi T$.
\item Si $A$ est un sous-objet de $X$, $S'$ n'est pas facteur de composition de $\Phi A$.
\item Si $B$ est un sous-quotient de $X$, $S'$ n'est pas facteur de composition de $\Phi B$.
\end{enumerate}
\end{lm}

\begin{proof} Il est trivial que (3) implique (1), et (2) implique (3) car $\Phi$ est exact à droite. Supposons maintenant (1) vérifié : une
récurrence sur la longueur montre que, pour tout sous-objet fini $F$ de $X$, $\Phi F$ n'a pas de facteur de composition $S'$ (utiliser encore l'exactitude à droite). Il suffit de passer à la
colimite pour conclure.\end{proof}

\begin{pr}
\label{prtech3}
Soient $X$ un objet localement fini de $\mathcal{A}$, $Y$ un sous-objet de $X$, $\Phi X\xrightarrow{\pi}B$ et $\Phi Y\xrightarrow{\pi'}B'$ des morphismes de $\mathcal{B}$ vérifiant les conditions suivantes :
\begin{enumerate}
\item le diagramme suivant commute :
$$\xymatrix{\Phi Y\ar[r]^{\pi'}\ar[d]_{\Phi i} & B'\ar[d]^{u} \\
\Phi X\ar[r]^{\pi} & B}$$
où $Y\xrightarrow{i} X$ désigne l'inclusion ;
\item l'objet $ker\,u$ n'a pas de facteur de composition $S'$ ;
\item le simple $S$ est $\Phi$-détecté par $S'$ dans $Y$ relativement
  à $\pi'$ ;
\item si $T$ est un facteur de composition de $X/Y$, alors $S'$ n'est pas facteur de composition de $\Phi T$.
\end{enumerate}
Alors $S$ est $\Phi$-détecté par $S'$ dans $X$ relativement à $\pi$.
\end{pr}

\begin{proof} Puisque $S'$ est facteur de composition de $im\,\pi'$ (par (3)), la condition (2) montre que $S'$ est facteur de composition de $im\,(u\circ\pi')$, c'est donc aussi le cas pour $im\,\pi$ (qui contient $im\,(u\circ\pi')$ par la condition (1)). 

Soit maintenant $A$ un sous-objet de $X$ tel que $S'$ est facteur de
composition de $\pi(im\,\Phi i)$ ; posons $A'=A\cap Y$. Comme $A/A'\hookrightarrow X/Y$, le lemme \ref{p-tech} prouve
que $\Phi(A/A')$ n'a pas de facteur de composition $S'$.

Considérons le diagramme commutatif
$$\xymatrix{\Phi A'\ar[r] \ar[d] & \Phi Y\ar[d] \ar[r]^{\pi'}& B'\ar[d]^{u}\\
\Phi A\ar[r] & \Phi X\ar[r]^{\pi} & B}$$
dont le carré de gauche est induit par les inclusions.
L'objet simple $S'$ est facteur de composition de $im\,(\Phi A'\to B)$ car 
$im\,(\Phi A\to B)/im\,(\Phi A'\to B)$ est un quotient de
$coker\,\Phi(A'\hookrightarrow A)\simeq\Phi(A/A')$, qui n'a pas de facteur $S'$. A fortiori, $S'$ est facteur de composition de
$\pi'\big(im\,(\Phi A'\to\Phi Y)\big)$ ; il s'ensuit (par (3)) que $S$ est facteur de composition de $A'$ (appliquer l'hypothèse (1)), donc de $A$, ce qui achève la démonstration. \end{proof}

\begin{pr}
\label{ptech2}
Soient $X$ un objet de $\mathcal{A}$,  $Y$ un sous-objet de $X$,
$\Phi X\xrightarrow{\pi}B$ et $\Phi(X/Y)\xrightarrow{\rho}C$ et $B\xrightarrow{v}
C$ des flèches de $\mathcal{B}$  vérifiant les conditions suivantes.
\begin{enumerate}
\item Le diagramme suivant commute.
\begin{equation}
\label{devv}\xymatrix{\Phi X \ar@{>>}[d] \ar[r]^\pi & B\ar[d]^v \\
\Phi(X/Y)\ar[r]^(.6){\rho} & C
}
\end{equation}
\item Le simple $S$ est $\Phi$-détecté par $S'$ dans $X/Y$
  relativement à $\rho$.
\item L'objet $ker\,v$ n'a pas de facteur de composition $S'$.
\end{enumerate}
Alors $S$ est $\Phi$-détecté par $S'$ dans $X$ relativement à $\pi$.
\end{pr}

\begin{proof} Le simple $S'$ est facteur de composition de $\pi(\Phi X)$, puisque cet objet
se projette sur $\rho(\Phi(X/Y))$. 

Soit à présent $A$ un sous-objet de $X$ tel que $S'$ est facteur de
composition de $\pi\big(im\,\Phi(A\hookrightarrow X)\big)$ ; on pose
$A'=A\cap Y$. L'examen du diagramme commutatif aux lignes exactes
$$\xymatrix{\Phi A'\ar[r]\ar[d] & \Phi A\ar@{>>}[r]\ar[d] &
  \Phi(A/A')\ar[d] \\
\Phi Y\ar[r] & \Phi X\ar@{>>}[r] & \Phi(X/Y)
}$$
montre que $im\,\big(\Phi(A/A')\to\Phi(X/Y)\big)$ est l'image de
$im\,(\Phi A\to\Phi X)$ par la projection $\Phi
X\twoheadrightarrow\Phi(X/Y)$. Par conséquent,
$\rho\big(im\,(\Phi(A/A')\to\Phi(X/Y))\big)$ est l'image de $\pi\big(im\,(\Phi
A\to\Phi X)\big)$ par $v$, et la dernière condition montre alors que $S'$
est facteur de composition de $\rho\big(im\,(\Phi(A/A')\to\Phi(X/Y))\big)$, de sorte que $A/A'$, et a fortiori $A$, a un facteur de
composition $S$, ce qu'il fallait démontrer.\end{proof}

\begin{cor}
\label{dw}
Soient $X$ un objet de $\mathcal{A}$,  $Y$ un sous-objet de $X$, et $\Phi X\xrightarrow{\pi}B$ un morphisme de $\mathcal{B}$ tels que :
\begin{enumerate}
\item $S'$ est facteur de composition de $im\,\pi$ ;
\item $S'$ n'est pas facteur de composition de $\Phi Y$ ;
\item $X/Y\simeq S$. 
\end{enumerate}
Alors $S$ est $\Phi$-détecté par $S'$ dans $X$ relativement à $\pi$.
\end{cor}

\begin{proof} Notons $i$ l'inclusion $Y\hookrightarrow X$ et $v$ la projection $B\twoheadrightarrow C=B/im\,(\pi\circ\Phi i)$, de sorte que $\pi$ induit
un morphisme $\rho : \Phi(X/Y)\to C$ rendant commutatif le diagramme (\ref{devv}). Les deux premières hypothèses montrent que $ker\,v$ n'a pas de facteur de composition
$S'$, tandis que $im\,\rho$ en a un, donc que $S$ est $\Phi$-détecté dans $X/Y$ relativement à $\rho$ grâce à la dernière hypothèse. La conclusion découle donc de la proposition \ref{ptech2}. 
\end{proof}

\subsection{Détection de facteurs de degré maximal dans un foncteur
  fini}\label{dfgm}

Comme le foncteur $(-:\Lambda^1)$ est un adjoint à gauche, il commute aux
colimites, ce qui permet de lui appliquer les résultats du paragraphe \ref{sectf}.

\begin{defi}\label{dfddf}
Soient $\lambda$ une partition régulière de longueur $r>0$ et
$X\in {\rm Ob}\,\F$. Nous dirons que $\lambda$ est
\textbf{$\Lambda^1$-détectable dans
  $X$} si, selon la terminologie de la définition \ref{dedet}, $S_{\lambda}$ est
$(-:\Lambda^1)$-détecté par $S_{\lambda^-_r}$ dans $X$.

Si $\pi : (X:\Lambda^1)\to B$ est un morphisme de $\F$,
nous dirons que $\lambda$ est \textbf{$\Lambda^1$-détectable dans $X$
  relativement à $\pi$} si $S_{\lambda}$ est
$(-:\Lambda^1)$-détecté par $S_{\lambda^-_r}$ dans $X$
relativement à $\pi$.
\end{defi}

\begin{nota} Soit $\lambda$ une partition de longueur $r>0$. On désigne par $\varsigma_{\lambda}: (W_{\lambda}:\Lambda^1)\twoheadrightarrow W_{\lambda^-_r}$
la projection donnée par la proposition \ref{cbdiv}.
\end{nota}

\begin{lm}\label{liodiot}
Soient $\lambda$ et $\mu$ deux partitions régulières, et
$r=l(\lambda)$ ; on suppose $r>0$. Si $\mu\vdash\Lambda^\lambda$ et $\lambda^-_r\vdash (S_\mu:\Lambda^1)$, alors $|\lambda|=|\mu|$ et $\mu\leq\lambda^{+,-}_{1,r}$.
\end{lm}

\begin{proof} On a $|\lambda|\geq|\mu|$ car $\mu\vdash \Lambda^{\lambda}$, et
$|\lambda|-1\leq|\mu|-1$ car $\lambda^-_r\vdash (S_{\mu}:\Lambda^1)$, d'où $|\lambda|=|\mu|$. La proposition \ref{fcsdl} entraîne maintenant $\mu^-_{1}\leq\lambda^-_r$, d'où $\mu\leq\lambda^{+,-}_{1,r}$.\end{proof}

\begin{pr}
\label{dew}
Soit $\lambda$ une partition régulière Weyl-séparante non nulle. Alors $\lambda$ est $\Lambda^1$-détectable dans $W_{\lambda}$ relativement à $\varsigma_{\lambda}$.
\end{pr}

\begin{proof} On applique le corollaire \ref{dw} avec $X=W_{\lambda}$
  et $Y={\rm rad}\,W_{\lambda}$, de sorte que la dernière condition
  est satisfaite grâce au théorème/définition \ref{repf}. La première est vérifiée parce que
$\lambda^-_r\vdash W_{\lambda^-_r}=im\,\varsigma_\lambda$.

La seconde condition du  corollaire \ref{dw} provient de l'hypothèse de Weyl-séparation, via le lemme \ref{p-tech}. En effet, supposons qu'elle ne soit pas satisfaite : il existerait une partition régulière
$\mu$ telle que $\mu\vdash {\rm rad}\, W_{\lambda}$ et
$\lambda^-_r\vdash (S_{\mu}:\Lambda^1)$, d'où $\mu\vdash {\rm rad}\, W_{\lambda}$ et $\mu\leq\lambda^{+,-}_{1,r}$ par le lemme \ref{liodiot}, en contradiction avec le fait que $\lambda$ est W-séparante.
\end{proof}

\begin{cor}\label{nvco}
On conserve les hypothèses de la proposition \ref{dew}. Soit $X$ un sous-objet
de $\Lambda^\lambda$ tel que :
\begin{itemize}
\item $W_\lambda\subset X$,
\item il n'existe pas de partition régulière $\mu$ de $|\lambda|$ telle que $\mu\vdash X/W_\lambda$ et $\mu\leq\lambda^{+,-}_{1,r}$.
\end{itemize}
Alors $\lambda$ est $\Lambda^1$-détectable dans $W_{\lambda}$ relativement à la composée $$\pi_X : (X:\Lambda^1)\to (\Lambda^{\lambda}:\Lambda^1)\twoheadrightarrow\Lambda^{\lambda^-_r}.$$
\end{cor}

\begin{proof} On applique la proposition \ref{prtech3} au sous-objet $W_\lambda$ de $X$. Ses deux premières conditions sont satisfaites, car le diagramme
$$\xymatrix{(W_\lambda:\Lambda^1)\ar[d]\ar@{>>}[r]^-{\varsigma_\lambda} & W_{\lambda^-_r}\ar@{^{(}->}[d] \\
(X:\Lambda^1)\ar[r]^-{\pi_X} & \Lambda^{\lambda^-_r}
}$$
commute.  

La troisième hypothèse de la proposition \ref{prtech3} est vérifiée par la proposition \ref{dew}, la dernière par le lemme \ref{liodiot}.
\end{proof}

\subsection{Détection dans $\bar{I}^{\otimes r}\otimes\Lambda^n$}

En vue d'appliquer les résultats précédents à la détection de sous-foncteurs de $\bar{I}^{\otimes r}\otimes \Lambda^\lambda$, nous établissons deux lemmes simples qui permettront de passer de la détection dans une partie homogène (à laquelle la section précédente est adaptée) à la détection \guillemotleft\, globale \guillemotright. Rappelons que l'on a un isomorphisme
$$p^{hom}_m(\bar{I}^{\otimes r})\simeq\underset{\alpha_1+\dots+\alpha_r=m}{\bigoplus_{\alpha_1,\dots,\alpha_r>0}}\Lambda^\alpha$$
pour tout entier $m$, via lequel nous identifierons souvent les deux membres.

\begin{lm}
\label{techdeg}
Soient $r$, $n$, $k$ trois entiers strictement positifs, et $\alpha$,
$\beta$ deux partitions régulières telles que 
$\alpha\vdash p^{hom}_k(\bar{I}^{\otimes r}\otimes\Lambda^n)$, $\beta\vdash
(S_\alpha : \Lambda^1)$ et $|\beta|<k-1$. Alors $l(\beta)\leq r$ et $\beta_{r}\leq n$.  
\end{lm}

\begin{proof} Supposons d'abord $|\alpha|<k$. On a  alors $l(\alpha)<r+1$ et $\alpha_{r}\leq n$ par le théorème \ref{fclambda}. En utilisant la proposition \ref{fcsdl}, on obtient $l(\beta)\leq l(\alpha)\leq r$ et $\beta_{r}\leq\alpha_r\leq n$.

Supposons désormais $|\alpha|=k$ : on a donc $|\beta|<|\alpha|-1$, et
la proposition \ref{fcsdl} donne $l(\beta)<l(\alpha)$ et
$\beta_{r}\leq\alpha_{r+1}$, donc, par le théorème \ref{fclambda},
on a $l(\beta)\leq r$ et $\beta_{r}\leq n$.\end{proof}

\begin{lm}
\label{homdet}
Soient $r$, $m$, $k$ des entiers strictement positifs,  $X$ un
sous-objet de $\bar{I}^{\otimes r}\otimes\Lambda^k$ et $\mu$ une partition régulière de $m$ telle que $\mu_{r+1}=k$. 
Si $\mu$ est $\Lambda^1$-détectable dans $p^{hom}_{m}(X)$ relativement au morphisme
$$(p_m^{hom}(X):\Lambda^1)\to (p_m^{hom}(\bar{I}^{\otimes r}\otimes\Lambda^k):\Lambda^1)\twoheadrightarrow p_{m-k}^{hom}(\bar{I}^{\otimes r})\otimes \Lambda^{k-1}$$
dont la première flèche est induite par l'inclusion et la seconde de
la proposition \ref{deriv} et de l'isomorphisme
$p_m^{hom}(\bar{I}^{\otimes r}\otimes\Lambda^k)\simeq p_{m-k}^{hom}(\bar{I}^{\otimes r})\otimes \Lambda^k$, alors $\mu$ est
$\Lambda^1$-détectable dans $X$ relativement au morphisme 
$$(X:\Lambda^1)\to (\bar{I}^{\otimes
  r}\otimes\Lambda^k:\Lambda^1)\simeq\bar{I}^{\otimes
  r}\otimes\Lambda^{k-1}$$ induit par l'inclusion via le corollaire \ref{cr-comdl}.
\end{lm}

\begin{proof} La proposition \ref{ptech2} prouve que $\mu$ est $\Lambda^1$-détectable dans $p_{m}X$ relativement au morphisme $(p_m X:\Lambda^1)\to (p_m(\bar{I}^{\otimes r}\otimes\Lambda^k):\Lambda^1)\twoheadrightarrow p_{m-k}(\bar{I}^{\otimes r})\otimes\Lambda^{k-1}$. En effet, le diagramme 
$$\xymatrix{(p_m X:\Lambda^1)\ar@{>>}[d]\ar[r] & (p_m(\bar{I}^{\otimes r}\otimes\Lambda^k):\Lambda^1)\ar@{>>}[r]\ar@{>>}[d] & p_{m-k}(\bar{I}^{\otimes r})\otimes \Lambda^{k-1}\ar@{>>}[d] \\
(p_m^{hom}(X):\Lambda^1)\ar[r]  & (p_m^{hom}(\bar{I}^{\otimes r}\otimes\Lambda^k):\Lambda^1)\ar@{>>}[r] & p_{m-k}^{hom}(\bar{I}^{\otimes r})\otimes\Lambda^{k-1}
}$$
commute, et le noyau de la flèche verticale de droite est de degré $<m-1$, donc sans facteur de composition $S_{\mu^-_{r+1}}$.

On termine la démonstration en utilisant la proposition \ref{prtech3} avec $Y=p_{m}X$ : le diagramme
$$\xymatrix{(p_m X:\Lambda^1)\ar[d]\ar[r] & (p_m(\bar{I}^{\otimes r}\otimes\Lambda^k):\Lambda^1)\ar@{>>}[r]\ar[d] & p_{m-k}(\bar{I}^{\otimes r})\otimes\Lambda^{k-1} \ar[d] \\
(X:\Lambda^1)\ar[r]  & (\bar{I}^{\otimes r}\otimes\Lambda^k:\Lambda^1)\ar[r]^-\simeq & \bar{I}^{\otimes r}\otimes\Lambda^{k-1}
}$$
commute (où les flèches verticales sont induites par les inclusions),
et la flèche verticale de droite est injective, ce qui montre que les deux
premières hypothèses de ladite proposition sont vérifiées. Nous venons
de voir que la troisième l'est ; quant à la dernière, elle provient du
lemme précédent : si elle était en défaut, on disposerait de $i>m$
d'une partition régulière $\alpha$ telle que $\alpha\vdash
p^{hom}_i(\bar{I}^{\otimes r}\otimes\Lambda^k)$ et $\mu^-_{r+1}\vdash
(S_\alpha: \Lambda^1)$. Donc 
$(\mu^-_{r+1})_{r}\leq k$ par le lemme \ref{techdeg}. Mais par hypothèse $\mu_{r+1}=k$, d'où  $(\mu^-_{r+1})_{r}=\mu_{r}>k$, contradiction qui achève la démonstration.\end{proof}

\begin{pr}
\label{cr-fd1}
Soient $r$, $m$, $k$ des entiers strictement positifs,
$A$ un sous-foncteur de $\bar{I}^{\otimes r}$, $X$ un sous-foncteur de
$A\otimes\Lambda^k$, $\mu$ une partition régulière de longueur $r$ de
$m-k$ telle que la suite d'entiers $\lambda=(\mu_1,\dots,\mu_r,k)$ est
une partition régulière Weyl-séparante et $\alpha :
p^{hom}_m(X)\to\Lambda^{\lambda}$ un morphisme vérifiant les
propriétés suivantes.
\begin{enumerate}
\item Il existe un morphisme $\beta:p^{hom}_{m-k}(A)\to\Lambda^\mu$ tel que $\alpha$ coïncide avec la composée
$$p^{hom}_m(X)\hookrightarrow p^{hom}_m(A\otimes\Lambda^k)=p^{hom}_{m-k}(A)\otimes\Lambda^k\xrightarrow{\beta\otimes\Lambda^k}\Lambda^\mu\otimes\Lambda^k\,.$$
\item L'image de $\alpha$ contient $W_{\lambda}$.
\item Il n'existe pas de partition régulière $\nu$ de $m$ telle que $\nu\vdash im\,\alpha/W_{\lambda}$ et $\nu\leq\lambda^{+,-}_{1,r+1}$.
\item \label{hhh} Si $\nu$ est une partition régulière de $m-k$ telle que
  $\nu\vdash ker\,\beta$, alors $\nu>\mu$.
\end{enumerate}

\smallskip

Alors $\lambda$ est $\Lambda^1$-détectable dans $X$ relativement au morphisme
$$(X :\Lambda^1)\to(\bar{I}^{\otimes r}\otimes\Lambda^k:\Lambda^1)\simeq\bar{I}^{\otimes r}\otimes\Lambda^{k-1}.$$
\end{pr}

\begin{proof}  Le corollaire \ref{nvco} montre que $\lambda$ est
  $\Lambda^1$-détectable dans $im\,\alpha$ relativement au morphisme
  $(im\,\alpha:\Lambda^1)\to \Lambda^{k-1}$ induit par l'inclusion.
On prouve maintenant que $\lambda$ est $\Lambda^1$-détectable dans $p^{hom}_m(X)$ relativement au morphisme
$$(p^{hom}_m(X):\Lambda^1)\to (p^{hom}_m(A\otimes\Lambda^k):\Lambda^1)\twoheadrightarrow p^{hom}_{m-k}(A)\otimes\Lambda^{k-1}$$
en employant la proposition \ref{ptech2}, avec le sous-objet
$ker\,\alpha$ et $u=\beta\otimes\Lambda^{k-1}$.

Pour la première condition, on constate que le diagramme 
$$\xymatrix{(p^{hom}_m(X):\Lambda^1)\ar[r]\ar@{>>}[d]_{(\alpha:\Lambda^1)}
  &
  (p^{hom}_m(A\otimes\Lambda^\lambda):\Lambda^1)\ar@{>>}[r]\ar[d]_{(\beta\otimes\Lambda^\lambda:\Lambda^1)} &  p^{hom}_{m-k}(A)\otimes\Lambda^{k-1}\ar[d]^-{\beta\otimes\Lambda^{k-1}} \\
(im\,\alpha:\Lambda^1)\ar[r] & (\Lambda^\mu\otimes\Lambda^\lambda : \Lambda^1)\ar@{>>}[r] & \Lambda^\mu\otimes \Lambda^{k-1}
}$$
commute.

Nous avons montré précédemment que la deuxième hypothèse de ladite proposition est vérifiée.

Pour la dernière, il s'agit d'établir que $S_{\lambda^-_{r+1}}$
n'est pas facteur de composition de
$ker\,\beta\otimes\Lambda^{k-1}$. C'est une conséquence
directe de l'hypothèse (\ref{hhh}) et du théorème \ref{fclambda}
: si $\nu'$ est une partition régulière de $m-1$ telle que
  $\nu'\vdash ker\,\beta\otimes\Lambda^{k-1}$, alors
  $\nu'\vdash\Lambda^{(\nu,k-1)}$, donc
  $\nu'\geq (\nu,k-1)$, où $\nu$ est une partition régulière de $m-k$ telle que
  $\nu\vdash ker\,\beta$.

Par conséquent (cf. remarque \ref{remf}), $\lambda$ est $\Lambda^1$-détectable dans $p^{hom}_m(X)$ relativement au morphisme
$$(p^{hom}_m(X):\Lambda^1)\to (p^{hom}_m(\bar{I}^{\otimes r}\otimes\Lambda^k):\Lambda^1)\twoheadrightarrow p^{hom}_{m-k}(\bar{I}^{\otimes r})\otimes\Lambda^{k-1}.$$ 

\smallskip

La conclusion résulte maintenant du lemme \ref{homdet}.\end{proof}

\section{Application à la structure de $I^{\otimes
    2}\otimes\Lambda^n$}\label{sec4}

Comme pour les avancées déjà connues dans l'étude de la conjecture
artinienne (cf. \cite{Piriou} et \cite{Po3}), la stratégie de la détermination de la structure des foncteurs
$I^{\otimes 2}\otimes\Lambda^n$, incluant le théorème \ref{thp}, que
nous mettons en \oe uvre dans cette section consiste en deux pas :
\begin{enumerate}\item réduire l'étude de ces foncteurs à celle de
  foncteurs plus simples ;
\item montrer que ces derniers n'ont pas de sous-foncteur propre \go
  trop gros\gf.
\end{enumerate}
La première étape sera réalisée par des constructions explicites liées
aux représentations des groupes symétriques, tandis que la seconde
repose sur l'utilisation du foncteur de division par $\Lambda^1$.

Afin de \go dévisser\gf au maximum les foncteurs
$\bar{I}^{\otimes 2}\otimes\Lambda^n$, on commence par ramener l'étude
de $\bar{I}^{\otimes 2}$ à celle de foncteurs plus simples. La décomposition des injectifs standard en somme directe
d'injectifs indécomposables peut être raffinée efficacement à l'aide
des notions de {\em foncteur co-Weyl} et de {\em filtration
  J-bonne} introduites par Powell dans \cite{Po1}, dont nous n'aurons pas
explicitement usage. En effet, dans le cas de
$\bar{I}^{\otimes 2}$, on obtient très simplement une filtration
explicite. Nous rappelons, dans le premier paragraphe, ces
considérations, et donnons des propriétés des facteurs de composition
des \go briques élémentaires\gf (autres que le foncteur $\bar{I}$) de
$\bar{I}^{\otimes 2}$, les foncteurs $L(2)$ et $\bar{D}(2)$. Nous verrons ainsi
qu'ils sont \go engendrés\gf par des foncteurs simples associés à des
partitions {\em alternées}.
%qui apparaissent très explicitement. 

\begin{rem} La possibilité de généraliser ces résultats aux foncteurs co-Weyl
supérieurs pose rapidement des problèmes techniques assez ardus ; quelques
renseignements remarquables (mais peu explicites) sur les facteurs de composition de ces
foncteurs  sont toutefois donnés en toute généralité dans  \cite{these}.
\end{rem}

Le dévissage de $\bar{I}^{\otimes 2}\otimes\Lambda^n$ obtenu par
tensorisation par $\Lambda^n$ de celui de $\bar{I}^{\otimes 2}$
précédemment évoqué ne
s'avère pas suffisant pour comprendre la structure de ce foncteur. Pour étudier le foncteur
$\bar{I}\otimes\Lambda^n$, on en définit d'abord un \go bon\gf
sous-foncteur $\bar{K}_n$, qui est l'image de la flèche de but
$\bar{I}\otimes\Lambda^n$ de la suite exacte longue
$$\cdots\to\bar{I}\otimes\Lambda^n\to\bar{I}\otimes\Lambda^{n-1}\to\dots\to\bar{I}\otimes\Lambda^1\to\bar{I}\to
0$$
(cf. \cite{Piriou}). De même, des suites exactes permettent
d'introduire des sous-foncteurs adéquats $L^2_n$ et $D^2_n$ de $L(2)\otimes\Lambda^n$ et
$\bar{D}(2)\otimes\Lambda^n$ respectivement, qui apparaissent à la
fois comme
noyau et image de flèches explicites, ce qui permet d'établir
aisément les propriétés nécessaires sur leurs facteurs
de composition à partir de celles de $L(2)$ et $\bar{D}(2)$. Cela fait l'objet du deuxième
paragraphe.

Enfin, le dernier paragraphe applique la proposition \ref{cr-fd1} sur la
détection de facteurs de composition par division par $\Lambda^1$ pour
en déduire le caractère artinien de type $2$
de $I^{\otimes 2}\otimes\Lambda^n$, par un argument de récurrence dont
l'initialisation fournie par
l'article \cite{P-art2}.

\subsection{La décomposition $\Lambda^2(\bar{I})\simeq L(2)\oplus\bar{D}(2)$}

Il est plus agréable de décrire le scindement de $\Lambda^2(\bar{P})$ dual de celui indiqué par le titre
de ce paragraphe, que l'on obtient à partir du fait suivant, qui
résulte d'un calcul direct.

\begin{lm}
Le morphisme $\Pi : \Lambda^2(\bar{P})\to\Lambda^2(\bar{P})$ donné par
$$ [u]\wedge [v]\mapsto \big([u]+[v]\big)\wedge [u+v]\qquad (u,v\in
V\setminus\{0\} ;\, V\in {\rm Ob}\,\E^f)$$
est un projecteur.
\end{lm}

\begin{rem}\label{rci} Soit $U$ l'endomorphisme $A\mapsto (f\mapsto
  A(f\circ\iota))$ de
  $I_{\mathbb{F}_2^{\oplus 2}}$, où $\iota$ désigne l'endomorphisme $(a,b)\mapsto
  (a,a+b)$ de $\mathbb{F}_2^{\oplus 2}$. Il est dual de
  l'endomorphisme de $P^{\otimes 2}$ donné par $[u]\otimes [v]\mapsto
  [u]\otimes [u+v]$. Le projecteur $\Pi$ peut se voir comme la composée
$$\Lambda^2\bar{P}\hookrightarrow\bar{P}^{\otimes 2}\xrightarrow{DU}\bar{P}^{\otimes 2}\twoheadrightarrow\Lambda^2\bar{P}.$$
\end{rem}

\begin{defi}
On pose $P_{2,1}=im\,\Pi$, $\bar{G}(2)=ker\,\Pi$, $L(2)=DP_{2,1}$ et
$\bar{D}(2)=D\bar{G}(2)$. Ainsi $\Lambda^2(\bar{P})\simeq
P_{2,1}\oplus \bar{G}(2)$ et $\Lambda^2(\bar{I})\simeq L(2)\oplus\bar{D}(2)$.
\end{defi}

\begin{rem}\label{fcli} Ces foncteurs peuvent être caractérisés comme
  suit --- on pourra se
référer à \cite{P-art2}, §\,$1.1$ à ce sujet. 
\begin{itemize}
\item Le foncteur $P_{2,1}$ est la couverture projective de $S_{(2,1)}$.
\item Le foncteur $\bar{G}(2)$ est isomorphe au foncteur $\FF[Gr_2]$,
  où $Gr_2(V)$ désigne la grassmannienne des plans d'un espace vectoriel $V$, l'action sur les morphismes étant l'action donnée par $\FF[Gr_2](f)\big([\pi]\big)=[f(\pi)]$ si $f(\pi)$ est un plan, $\FF[Gr_2](f)\big([\pi]\big)=0$ sinon.

 Plus précisément, si $(u,v)$ est une famille libre de $V\in {\rm Ob}\,\E^f$, dénotons par $<u,v>\in Gr_2(V)$ le plan qu'elle engendre. Alors le morphisme $\Lambda^2(\bar{P})\xrightarrow{p}\FF[Gr_2]$ donné par
$[u]\wedge [v]\mapsto [<u,v>]$ se restreint en un isomorphisme de $\bar{G}(2)$ sur $\FF[Gr_2]$. 
\end{itemize}
\end{rem}

Avant d'indiquer les facteurs de composition que nous utiliserons pour la $\Lambda^1$-détection dans $L(2)\otimes\Lambda^n$ et $\bar{D}(2)\otimes\Lambda^n$, nous énonçons deux lemmes formels. 

\begin{lm}
Soient $\lambda$ une partition régulière et $X$ un foncteur
analytique. Supposons que $S_\lambda$ est facteur de composition unique de
$X$. Alors il existe un plus petit sous-objet $X[\lambda]$ de $X$ tel
que $\lambda\vdash X[\lambda]$. Le foncteur $X[\lambda]$ est fini, de
degré supérieur à $|\lambda|$. De plus, $S_\lambda$ est le cosocle de
$X[\lambda]$, et cette propriété caractérise $X[\lambda]$ parmi les
sous-objets finis de $X$.
\end{lm}

\begin{proof} Comme $X$ est analytique, $X$ a un
  sous-objet fini $F$ tel que $\lambda\vdash F$. D'autre part, si $A$ et $B$ sont deux sous-objets de $X$
  tels que $\lambda\vdash A$ et $\lambda\vdash B$, alors
  $\lambda\vdash A\cap B$, car sinon $\lambda$ serait facteur de
  composition (au moins) double de $(A\oplus B)/(A\cap B)\simeq
  A+B\subset X$. Par conséquent, l'intersection $X[\lambda]$ des sous-objets $A$ de $X$ tels que
  $\lambda\vdash A$ convient. 

On a $\deg\,X[\lambda]\geq \deg\,S_\lambda=|\lambda|$. Par ailleurs,
si $\pi : X[\lambda]\twoheadrightarrow S$ est un épimorphisme avec $S$
simple, $ker\,\pi$ n'a pas de facteur de composition
$S_\lambda$, donc $S\simeq S_\lambda$ ; on en déduit ${\rm cosoc}\,X[\lambda]=S_\lambda$. Réciproquement, si $Y$ est un
sous-objet de $X$ de cosocle $S_\lambda$, $X[\lambda]\subset Y$ ; si
l'inclusion était stricte, $S_\lambda$ serait facteur de composition
de ${\rm rad}\,Y\supset X[\lambda]$, donc serait facteur de composition (au moins) double de $Y$, contradiction qui achève la démonstration.\end{proof}

\begin{ex}[fondamental] Le foncteur simple $S_\lambda$ est facteur de composition unique dans $\Lambda^\lambda$ et $\Lambda^\lambda[\lambda]=W_\lambda$.
\end{ex}

\begin{defi}
Soient $X$ un foncteur analytique et $\lambda$ une partition régulière telle que $\lambda\vdash X$. Nous dirons que $S_\lambda$ est \textbf{bien placé} dans $X$ si $S_\lambda$ n'est pas facteur de composition de $X/p_{|\lambda|}(X)$. 
\end{defi}

\begin{rem} Si $S_\lambda$ est facteur de composition unique de $X$, cela équivaut à l'inclusion $X[\lambda]\subset p_{|\lambda|} X$, ou encore à $\lambda\vdash p^{hom}_{|\lambda|} X$.
\end{rem}

\begin{lm}
\label{detcor}
Soient $X$ un foncteur analytique et $\lambda$ une partition régulière
d'un entier $n$ telle que $S_\lambda$ est facteur de composition
unique bien placé de $X$. Supposons aussi que $Y\xrightarrow{f}X$ est
un morphisme de $\mathcal{F}_\omega$ et $A$ un sous-objet cohomogène
de degré $n$ de $p_n Y$ tels que $p^{hom}_n(f)(p^{hom}_n A)=(p^{hom}_n
X)[\lambda]$. On a alors $f(A)=X[\lambda]$.
\end{lm}

\begin{proof} Les foncteurs $p_i$ étant exacts à gauche et se plongeant naturellement dans le foncteur identité, on a $p^{hom}_n(f)(p^{hom}_n A)\twoheadrightarrow p^{hom}_n\big(f(A)\big)$. On en déduit $S_\lambda\simeq {\rm cosoc}\,p^{hom}_n(f)(p^{hom}_n A)\twoheadrightarrow{\rm cosoc}\,p^{hom}_n\big(f(A)\big)$.

D'autre part, comme $\deg A\leq n$, on dispose d'un épimorphisme
$f(A)\twoheadrightarrow p^{hom}_n(f)(p^{hom}_n A)$, donc aussi ${\rm cosoc}\,f(A)\twoheadrightarrow {\rm cosoc}\,p^{hom}_n(f)(p^{hom}_n A)$.

Enfin, le quotient $f(A)$ de $A$ est cohomogène de degré $n$ (ou nul), donc
la projection $f(A)\twoheadrightarrow p^{hom}_n\big(f(A)\big)$ induit
un isomorphisme ${\rm cosoc}\,f(A)\xrightarrow{\simeq} {\rm
  cosoc}\,p^{hom}_n\big(f(A)\big)$. Conséquemment, $S_\lambda\simeq {\rm cosoc}\,f(A)$, d'où le lemme. \end{proof}

Nous revenons aux foncteurs $\Lambda^2(\bar{I})$, $L(2)$ et
$\bar{D}(2)$, dont nous étudions les facteurs de composition à travers
leur filtration polynomiale.

Pour tout entier $n>0$, on a un isomorphisme
\begin{equation}
\label{fpol}
p^{hom}_{n}\big(\Lambda^2(\bar{I})\big)\simeq\left(\underset{a>b>0}{\bigoplus_{a+b=n}}\Lambda^{a}\otimes\Lambda^b\right)\oplus\Lambda^2(\Lambda^{n/2})
\end{equation}
où, par convention, le dernier terme est nul si $n$ est impair. Via
cette identification, le plongement de $p^{hom}_{n}\big(\Lambda^2(\bar{I})\big)$ dans
$$p^{hom}_{n}\big(\bar{I}^{\otimes 2}\big)\simeq\underset{a,b>0}{\bigoplus_{a+b=n}}(\Lambda^{a}\otimes\Lambda^b)$$
s'obtient comme somme des morphismes $\Lambda^a\otimes\Lambda^b\xrightarrow{id\oplus\tau}(\Lambda^a\otimes\Lambda^b)\oplus (\Lambda^b\otimes\Lambda^a)$, pour $a>b>0$ et $a+b=n$, $\tau$ désignant l'isomorphisme d'échange des deux facteurs du produit tensoriel, et de l'inclusion $\Lambda^2(\Lambda^{n/2})\hookrightarrow\Lambda^{n/2}\otimes\Lambda^{n/2}$.

\begin{lm}
\label{linterm}
Soient $i,j,k,l$ et $n$ des entiers vérifiant $i>j>0$, $k\geq l>0$ et
$i+j=k+l=n$. La composée $$\Lambda^i\otimes\Lambda^j\hookrightarrow
p^{hom}_{n}(\Lambda^2(\bar{I}))\xrightarrow{p^{hom}_{n}(D\Pi)}p^{hom}_{n}(\Lambda^2(\bar{I}))\twoheadrightarrow\Lambda^k\otimes\Lambda^l$$
dont les première et dernière flèches sont déduites de (\ref{fpol}) est la somme des
morphismes :
\begin{flushleft}
\begin{itemize}
\item $\theta_{i,j,k-i}$ si $k\geq i$ ;
\item $D\theta_{k,l,i-k}$ si $k\leq i$ ;
\item $\Lambda^i\otimes\Lambda^j\xrightarrow{\tau}\Lambda^j\otimes\Lambda^i\xrightarrow{\theta_{j,i,k-j}}\Lambda^k\otimes\Lambda^l$ si $k\geq j$, où $\tau$
désigne la flèche échangeant les deux facteurs du produit tensoriel ;
\item $\Lambda^i\otimes\Lambda^j\xrightarrow{\tau}\Lambda^j\otimes\Lambda^i\xrightarrow{D\theta_{k,l,j-k}}\Lambda^k\otimes\Lambda^l$ si $k\leq j$.
\end{itemize}
\end{flushleft}

En particulier, pour $(i,j)=(k,l)$, le morphisme en question est $\theta_{j,i,i-j}\circ\tau$.
\end{lm}

\begin{proof} Cela provient de la remarque \ref{rci}, en utilisant que $U$ induit au niveau de la
filtration polynomiale les morphismes $$\Lambda^i\otimes\Lambda^j\xrightarrow{\bigoplus_{0\leq t\leq j}\theta_{i,j,t}}\underset{0\leq t\leq j}{\bigoplus}
\Lambda^{i+t}\otimes\Lambda^{j-t}.$$  En effet, si
$a_1,\dots,a_i,b_1,\dots,b_j$ sont des éléments d'un espace vectoriel $V$, l'élément
$(a_1\wedge\dots\wedge a_i)\otimes (b_1\wedge\dots\wedge b_j)$ de $\Lambda^i(V)\otimes\Lambda^j(V)$ se relève en l'élément de
$p_{i+j}(I_{\mathbb{F}_2^{\oplus 2}})(V)$ donné par $(l,l')\mapsto\big(\prod_{r=1}^i l(a_r)\big)\big(\prod_{s=1}^j l'(b_s)\big)$ ($(l,l')\in (V^*)^2$ ; on
identifie $I_{\mathbb{F}_2^{\oplus 2}}(V)$ et $I^{\otimes 2}(V)$) ; ensuite développer le produit dans l'élément
$(l,l')\mapsto\big(\prod_{r=1}^i l(a_r)\big)\big(\prod_{s=1}^j (l(b_s)+l'(b_s))\big)$ de $p_{i+j}(I_{\mathbb{F}_2^{\oplus 2}})(V)$ qui est l'image du précédent
par le morphisme induit par $U$.
\end{proof}

\begin{lm}
\label{ldtcl} Soient $i>j>0$ des entiers. La restriction à $W_{(i,j)}$ du morphisme
$\Lambda^i\otimes\Lambda^j\hookrightarrow
p^{hom}_{i+j}\big(\Lambda^2(\bar{I})\big)\xrightarrow{p^{hom}_{i+j}(D\Pi)}p^{hom}_{i+j}\Lambda^2(\bar{I})\twoheadrightarrow\Lambda^{i}\otimes\Lambda^j$
coïncide avec l'identité.
\end{lm}

\begin{proof} Cela résulte des lemmes \ref{linterm} et \ref{mweyl}.\end{proof}

\begin{nota}\label{nfg} Soit $i\in\mathbb{N}^*$. On note $f_i$ la
  composée
$$p^{hom}_{2i+1}L(2)\hookrightarrow
  p^{hom}_{2i+1}\Lambda^2(\bar{I})\twoheadrightarrow\Lambda^{i+1}\otimes\Lambda^i\,,$$
et $g_i$ le morphisme
$$p^{hom}_{2i+3}\bar{D}(2)\hookrightarrow p^{hom}_{2i+3}\Lambda^2(\bar{I})\twoheadrightarrow\Lambda^{i+2}\otimes\Lambda^{i+1}\xrightarrow{\theta}\Lambda^{i+3}\otimes\Lambda^i\,.$$
\end{nota}

\begin{pr}
\label{detfcl} Soit $i\in\mathbb{N}^*$. 
\begin{enumerate}
\item \begin{enumerate}
\item Le foncteur $L(2)$ possède un unique facteur de composition $S_{(i+1,i)}$, tandis que $\bar{D}(2)$ n'en a pas.
\item De plus, celui-ci est bien placé. Précisément, $W_{(i+1,i)}\subset im\,f_i$. 
\item En revanche, $im\,f_i$ n'a pas de facteur de composition $S_{(i+2,i-1)}$.
\end{enumerate}
\item \begin{enumerate}
\item Le foncteur $\bar{D}(2)$ possède un unique facteur de composition $S_{(i+3,i)}$.
\item De plus, celui-ci est bien placé, et 
  $W_{(i+3,i)}\subset im\,g_i$. Par conséquent, $\Lambda^{i+2,i+1}[i+3,i]\subset p^{hom}_{2i+3}\bar{D}(2)$.
\item En revanche, $im\,g_i$ n'a pas de facteur de composition $S_{(i+4,i-1)}$.
\end{enumerate}
\end{enumerate}
\end{pr}

\begin{proof} Tout d'abord, (\ref{fpol}) montre que les facteurs de composition de $\bar{I}^{\otimes 2}$ sont tous bien placés.

D'autre part, pour $i\in\mathbb{N}^*$, $S_{(i+1,i)}$ est facteur de
composition unique de $p^{hom}_{2i+1}\Lambda^2(\bar{I})$, de sorte que
(1) (a) et  (1) (b) découlent du lemme \ref{ldtcl}.

\smallskip

Les partitions $(i+2,i+1)$ étant alternées, donc
Weyl-séparantes (proposition \ref{pfc}), $S_{(i+3,i)}$ n'est pas
facteur de composition de $W_{(i+2,i+1)}$, donc est facteur de
composition unique de $\Lambda^{i+2}\otimes\Lambda^{i+1}$, et
$\theta_{i+2,i+1,1}\big((\Lambda^{i+2}\otimes\Lambda^{i+1})[i+3,i]\big)=W_{(i+3,i)}$
(utiliser la filtration de Weyl usuelle de
$\Lambda^{i+2}\otimes\Lambda^{i+1}$ --- cf. \cite{Piriou}, §\,1.1). Par conséquent,
$\Lambda^2(\bar{I})$ a exactement deux facteurs de composition
$S_{(i+3,i)}$ pour $i\in\mathbb{N}^*$ ; le lemme \ref{ldtcl} (qui
montre notamment que l'un d'entre eux apparaît dans $L(2)$) montre
qu'il suffit de voir pour démontrer (2) (a) et  (2) (b) que $L(2)$ a un seul facteur de composition $S_{(i+3,i)}$.

Pour cela, on utilise l'endofoncteur $\tilde{\nabla}_2$ de $\mathcal{F}$ introduit par G. Powell, dont on emploiera les propriétés suivantes (pour la définition de $\tilde{\nabla}_2$ et la démonstration de ces propriétés, voir \cite{Po2}) :
\begin{enumerate}
\item $\tilde{\nabla}_2$ préserve les injections et les surjections. En particulier, si $\lambda\vdash X$, alors $\tilde{\nabla}_2 S_\lambda$ est un sous-quotient de $\tilde{\nabla}_2 X$.
\item $\tilde{\nabla}_2 S_{(i+1,j+1)}=S_{(i,j)}$ pour $i>j\geq 0$.
\item $\tilde{\nabla}_2 L(2)=L(2)\oplus\bar{I}$. %et $\tilde{\nabla}_2 \bar{D}(2)=\bar{D}(2)$.
\item Pour $i>0$ et $j>t\geq 0$, $\tilde{\nabla}_2$ transforme $\Lambda^i\otimes\Lambda^j\xrightarrow{\theta_{i,j,t}}\Lambda^{i+t}\otimes\Lambda^{j-t}$ en
$\Lambda^{i-1}\otimes\Lambda^{j-1}\xrightarrow{\theta_{i-1,j-1,t}}\Lambda^{i+t-1}\otimes\Lambda^{j-t-1}$.     
\end{enumerate}

Si pour un $i\in\mathbb{N}^*$, $L(2)$ avait deux facteurs de
composition $S_{(i+3,i)}$, on en déduirait que $L(2)$ a un facteur de
composition $\Lambda^3$, ce qui n'est pas le cas puisque l'unique
facteur de composition $\Lambda^3$ de $\Lambda^2(\bar{I})$ apparaît
dans
$p^{hom}_3\big(\Lambda^2(\bar{I})\big)=\Lambda^2\otimes\Lambda^1\simeq
S_{(2,1)}\oplus\Lambda^3$, or dans cette décomposition $p^{hom}_3
L(2)=S_{(2,1)}$ et $p^{hom}_3\bar{D}(2)=\Lambda^3$ (utiliser le lemme
\ref{linterm}). Cela établit (2) (a) et  (2) (b).

Cela montre également (1) (c) pour $i=1$. Le cas général s'en déduit encore via l'utilisation de $\tilde{\nabla}_2$, car d'après les propriétés rappelées ci-avant, ce foncteur transforme la flèche $f_{i+1}$ en $f_i$, donc $im\,f_{i+1}$ en $im\,f_i$, de sorte que si  $(i+3,i)\vdash im\,f_{i+1}$,  alors $(i+2,i-1)\vdash im\,f_{i}$.

De même, il suffit de démontrer (2) (c) pour $i=1$. Or le lemme \ref{ldtcl} montre que $g_1$ a la même image que $\Pi_{4,1}$ (cf. notation \ref{nfd} et proposition \ref{pproj}), or ce dernier s'identifie à la projection $\Lambda^4\otimes\Lambda^1\twoheadrightarrow S_{(4,1)}$, ce qui achève la démonstration.
\end{proof}

La définition que nous rappelons ci-après correspond est adaptée pour
préciser la structure des injectifs co-tf de $\F$. Ainsi, une version
forte de la conjecture artinienne postule que le foncteur $I^{\otimes
  n}$ est artinien de type $n$ pour tout $n\in\mathbb{N}$ --- cf. \cite{Po1}.

\begin{defi}[cf. \cite{Po1} et \cite{GP-gal}] On définit par
  récurrence sur $n\in\mathbb{N}$ la notion de foncteur {\em simple
    artinien de type $n$} (resp. {\em foncteur artinien de type $n$}). Un foncteur simple artinien de type $0$ est un foncteur simple ; un foncteur simple artinien de type $n+1$ est un foncteur qui n'est pas artinien de type $n$ mais dont tous les sous-objets stricts sont artiniens de type $n$. Un foncteur artinien de type $n$ est un foncteur qui possède une filtration finie dont les sous-quotients sont simples artiniens de type $\leq n$.
\end{defi}

On vérifie par récurrence sur $n$ qu'un objet artinien de type $n$ est
artinien. Nos arguments ultérieurs se fonderont sur le résultat suivant.

\begin{theo}[Powell]\label{lmgp}
Il existe une filtration croissante $(F_n)_{n\in\mathbb{N}^*}$ sur $L(2)$ (resp. $\bar{D}(2)$) telle que : 
\begin{itemize}
\item tout sous-objet strict de $L(2)$ (resp. $\bar{D}(2)$) est inclus
  dans l'un des $F_n$ ;
\item pour $n\in\mathbb{N}^*$ et un foncteur fini $A$,
  $A\otimes F_n$ est artinien de type $1$ ;
\item si $\lambda$ est une partition régulière telle que $\lambda\vdash F_n$, alors $l(\lambda)\leq 2$ et $\lambda_2\leq n$.
\end{itemize}

Par suite, les foncteurs $L(2)$ et $\bar{D}(2)$ sont simples artiniens de type $2$.
\end{theo}

La démonstration est fournie dans \cite{P-art2}, à combiner aux résultats
de \cite{Po3} pour le deuxième point.

\begin{cor}
\label{prgp}
\begin{enumerate}
\item Le foncteur $L(2)$ est la réunion filtrante sur $i\in\mathbb{N}^*$ des $L(2)[i+1,i]$.
\item Le foncteur $\bar{D}(2)$ est la réunion filtrante sur $i\in\mathbb{N}^*$ des $\bar{D}(2)[i+3,i]$.
\end{enumerate}
\end{cor}

\begin{proof} Il s'agit d'une conséquence immédiate de la proposition
  \ref{detfcl} et du théorème \ref{lmgp}. \end{proof}

\begin{rem} La première assertion implique qu'un sous-foncteur $F$ de $L(2)$ tel que
  $(i+1,i)\vdash F$ pour une infinité d'entiers $i$ est égal à $L(2)$.
\end{rem}

\subsection{Les foncteurs  $L_n^2$ et $D_n^2$}
%\subsection{Les sous-objets  $L_n^2$ et $D_n^2$ de $\Lambda^2(\bar{I})\otimes\Lambda^n$}

  Nous sommes en mesure de donner une description explicite des \go briques
élémentaires\gf des foncteurs
$\Lambda^2(\bar{I})\otimes\Lambda^k$. C'est la suite exacte
(\ref{segen}) qui permettra de dévisser itérativement ces foncteurs à
l'aide des foncteurs $L^2_i$ et $D^2_i$ que nous allons définir et étudier.

\begin{conv} Dans ce paragraphe, $n$ désigne un entier strictement positif.
\end{conv}

Les constructions que nous allons exposer reposent sur la considération des
morphismes suivants. Là
encore, il est commode de commencer par introduire les flèches duales
de celles qui nous intéressent.

\begin{nota} On désigne par $g_n : \Lambda^2(\bar{P})\otimes\Lambda^n\to\Lambda^2(\bar{P})\otimes\Lambda^n(\Lambda^3)$ le morphisme donné par
$$\big([u]\wedge [v]\big)\otimes (a_1\wedge\dots\wedge a_n)\mapsto \big([u]\wedge [v]\big)\otimes\bigwedge_{1\leq i\leq n}(u\wedge v\wedge a_i)$$
et par $h_n : \bar{P}^{\otimes 2}\otimes\Lambda^{n-1}\to\Lambda^2(\bar{P})\otimes\Lambda^n$ le morphisme donné par
$$\big([u]\otimes [v]\big)\otimes (a_1\wedge\dots\wedge
a_{n-1})\mapsto \big([u]\wedge [v]\big)\otimes (u\wedge
a_1\wedge\dots\wedge a_{n-1}).$$
\end{nota}

\begin{lm}\label{lmdg} Le diagramme 
$$\xymatrix{\Lambda^2(\bar{P})\otimes\Lambda^n\ar[r]^-{g_n}\ar[d]_{\Pi\otimes\Lambda^n} & \Lambda^2(\bar{P})\otimes\Lambda^n(\Lambda^3)\ar[d]^{\Pi\otimes\Lambda^n(\Lambda^3)} \\
\Lambda^2(\bar{P})\otimes\Lambda^n\ar[r]^-{g_n} & \Lambda^2(\bar{P})\otimes\Lambda^n(\Lambda^3)
}$$
commute. Ainsi, $g_n$ s'identifie, via la décomposition
$\Lambda^2(\bar{P})\simeq P_{2,1}\oplus \bar{G}(2)$, à la somme directe de deux morphismes $P_{2,1}\otimes\Lambda^n\xrightarrow{g_n^L} P_{2,1}\otimes\Lambda^n(\Lambda^3)$ et $\bar{G}(2)\otimes\Lambda^n\xrightarrow{g_n^D}\bar{G}(2)\otimes\Lambda^n(\Lambda^3)$.
\end{lm}

Ce résultat provient d'un calcul direct.

Le lemme suivant utilise les morphismes $w$ définis au
corollaire \ref{crweyl}.

\begin{lm}\label{jvlt} Soient $i>j>0$ des entiers. La restriction à \\
$\Lambda^i\otimes\Lambda^j\otimes\Lambda^n(\Lambda^3)$ de $p^{hom}_{i+j+3n}(Dg_n)$ s'identifie à $w^{i,j,0}_n$ (composée avec l'inclusion $\Lambda^{i+n}\otimes\Lambda^{j+n}\otimes\Lambda^n\hookrightarrow p^{hom}_{i+j+3n}(\Lambda^2(\bar{I})\otimes\Lambda^n)$). En conséquence, $p^{hom}_{i+j+3n}(Dg_n)(W_{(i,j)}\otimes\Lambda^n(\Lambda^3))=W_{(i+n,j+n,n)}$.
\end{lm}

\begin{proof} La construction des morphismes $w^{i,j,0}_n$ et $g_n$
  montre qu'il suffit de traiter le cas $n=1$. La flèche $\Lambda^{i+1,j+1,1}\to\Lambda^{i,j,3}$
%$\Lambda^{i+1}\otimes\Lambda^{j+1}\otimes\Lambda^1\to\Lambda^i\otimes\Lambda^j\otimes\Lambda^3$
duale de $w^{i,j,0}_1=j_{(i,j),3}$ --- morphisme de la proposition \ref{prpra} --- s'obtient à partir du coproduit $\Lambda^{t+1}\to\Lambda^t\otimes\Lambda^1$ pour $t=i,j$ et du produit $\Lambda^1\otimes\Lambda^1\otimes\Lambda^1\to\Lambda^3$. La somme directe pour $i+j=n$ et $i,j>0$ de ces morphismes correspond donc à la filtration polynomiale de la flèche $\bar{P}^{\otimes 2}\otimes\Lambda^1\to\bar{P}^{\otimes 2}\otimes\Lambda^3$ donnée par $[u]\otimes [v]\otimes a\mapsto [u]\otimes [v]\otimes (u\wedge v\wedge a)$ ; par suite, la somme directe pour $i+j=n$ et $i>j>0$ de ces morphismes décrit la filtration polynomiale de $g_1 : \Lambda^2(\bar{P})\otimes\Lambda^n\to\Lambda^2(\bar{P})\otimes\Lambda^3$.
\end{proof}

Ce résultat nous permet de donner les principales propriétés
nécessaires à la détection de facteurs composition dans
$L(2)\otimes\Lambda^n$ et $\bar{D}\otimes\Lambda^n$.

\begin{pr}\label{ptdg}
\begin{enumerate}
\item \label{pot1} Pour tout entier $i>0$, $L(2)\otimes\Lambda^n$ contient un unique facteur de composition $S_{(i+n+1,i+n,n)}$. De plus,
$$(L(2)\otimes\Lambda^n)[i+n+1,i+n,n]=Dg_n^L\big(L(2)[i+1,i]\otimes\Lambda^n(\Lambda^3)\big)\,.$$
\item \label{pot2} Pour tout entier $i>0$, $\bar{D}(2)\otimes\Lambda^n$ contient un unique facteur de composition $S_{(i+n+3,i+n,n)}$. De plus,
$$(\bar{D}(2)\otimes\Lambda^n)[i+n+3,i+n,n]=Dg_n^D\big(\bar{D}(2)[i+3,i]\otimes\Lambda^n(\Lambda^3)\big)\,.$$
\item \label{segen2} La suite suivante est exacte. 
\begin{equation}\label{segen}
\bar{P}^{\otimes 2}\otimes\Lambda^{n-1}\xrightarrow{h_n}\Lambda^2(\bar{P})\otimes\Lambda^n\xrightarrow{g_n}\Lambda^2(\bar{P})\otimes\Lambda^n(\Lambda^3)
\end{equation}
\end{enumerate}
\end{pr}

\begin{proof} Le foncteur $L(2)[i+1,i]\otimes\Lambda^n(\Lambda^3)$ est
  cohomogène comme produit tensoriel de deux foncteurs cohomogènes
  (cf. corollaire \ref{crch}) de degré $3i+2n+1$, et sa partie
  homogène de degré $3n+2i+1$ est
  $W_{(i+1,i)}\otimes\Lambda^n(\Lambda^3)$ par la proposition
  \ref{detfcl}, de sorte que
$$p^{hom}_{3n+2i+1}\big(Dg_n^L\big)\big(p^{hom}_{3n+2i+1}(L(2)[i+1,i]\otimes\Lambda^n(\Lambda^3))\big)$$
$$=p^{hom}_{3n+2i+1}\big(Dg_n^L\big)(W_{(i+1,i)}\otimes\Lambda^n(\Lambda^3))=W_{(i+n+1,i+n,n)}$$
$$=p^{hom}_{3n+2i+1}\big(L(2)\otimes\Lambda^n\big)[i+n+1,i+n,n].$$ 
On a fait usage du lemme \ref{jvlt} pour la deuxième égalité. Le lemme
\ref{detcor} fournit donc l'assertion (\ref{pot1}).

On établit de même l'assertion (\ref{pot2}), en utilisant
également le dernier point du corollaire \ref{crweyl}.

L'assertion (\ref{segen2}) découle quant à elle d'un calcul direct.
\end{proof}

\begin{defi}
On pose $L^2_n=im\,Dg_n^L$ et $D^2_n=im\,Dg_n^D$. On note également $L^2_0=L(2)$ et $D^2_0=\bar{D}(2)$.
\end{defi}

\begin{cor}\label{cdetl}
\begin{itemize}
\item Le foncteur $L^2_n$ est la réunion filtrante sur
  $i\in\mathbb{N}^*$ des sous-foncteurs $(L(2)\otimes\Lambda^n)[i+n+1,i+n,i]$.
\item Le foncteur $D^2_n$ est la réunion filtrante sur $i\in\mathbb{N}^*$ des sous-foncteurs $(\bar{D}(2)\otimes\Lambda^n)[i+n+3,i+n,i]$.
\end{itemize}
\end{cor}

\begin{proof} C'est une conséquence immédiate du corollaire \ref{prgp}
  et de la proposition \ref{ptdg}.\end{proof}

Le lemme suivant est une variation sur le corollaire \ref{crweyl}
adapté au cas de $D^2_n$, légèrement plus technique que celui de
$L^2_n$, pour lequel ce corollaire suffira à nos investigations ultérieures.

\begin{lm}\label{plfcd}
Soient $k$ et $n$ des entiers strictements positifs. On a
$$w_n^{k+2,k+1,0}(\Lambda^{k+2,k+1}[k+3,k]\otimes\Lambda^n(\Lambda^3))=\Lambda^{k+n+2,k+n+1,n}[k+n+3,k+n,n].$$
\end{lm}

\begin{proof} Elle est entièrement analogue à celle de la proposition
  \ref{prpra}, en notant que pour tout $V\in {\rm Ob}\,\E^f$, $\Lambda^{k+2,k+1}[k+3,k](V)$ est le sous--espace vectoriel de $\Lambda^{k+2,k+1}(V)$ engendré par les éléments du type
$(a_1\wedge\dots\wedge a_k\wedge b\wedge c)\otimes (a_1\wedge\dots\wedge a_k\wedge d)+(a_1\wedge\dots\wedge a_k\wedge b\wedge d)\otimes (a_1\wedge\dots\wedge a_k\wedge c)+(a_1\wedge\dots\wedge a_k\wedge c\wedge d)\otimes (a_1\wedge\dots\wedge a_k\wedge b)$
pour $a_1,\dots,a_k,b,c,d\in V$.
\end{proof}

On rappelle que les morphismes $f_k$ et $g_k$ ont été introduits dans
la notation \ref{nfg}.

\begin{nota} Soient $n$ et $k$ des entiers tels que $k>n>0$. On pose
  $A_{k,n}=(f_k\otimes\Lambda^n)\big(p^{hom}_{2k+n+1}(L^2_n)\big)$ et $B_{k,n}=(g_k\otimes\Lambda^n)\big(p^{hom}_{2k+n+3}(D^2_n)\big)$.
\end{nota}

La proposition suivante contient tous les préliminaires nécessaires à
la détection de facteurs de composition dans les foncteurs $L^2_n$ et $D^2_n$.

\begin{pr}\label{lpcfc} Soient $n$ et $k$ des entiers tels que $k>n>0$.
\begin{enumerate}\item \begin{enumerate}\item\label{mri4} On a $W_{(k+1,k,n)}\subset A_{k,n}$.
\item\label{mri} Il n'existe pas de partition régulière $\nu$ de $2k+n+1$ telle
  que $\nu\vdash A_{k,n}/W_{(k+1,k,n)}$ et  $\nu\leq (k+2,k,n-1)$.
\end{enumerate}
\item \begin{enumerate}\item\label{mri5} On a $W_{(k+3,k,n)}\subset B_{k,n}$.
\item\label{mri7} Il n'existe pas de partition régulière $\nu$ de $2k+n+3$ telle
  que $\nu\vdash B_{k,n}/W_{(k+3,k,n)}$ et  $\nu\leq (k+4,k,n-1)$.
\end{enumerate}
\end{enumerate}
\end{pr}

\begin{proof} L'assertion (\ref{mri4}) s'obtient en
  combinant le lemme \ref{jvlt} et la proposition
  \ref{detfcl}. L'assertion (\ref{mri5}) s'établit pareillement,
  en utilisant aussi le lemme \ref{plfcd}. 

Montrons l'assertion (\ref{mri}). Grâce à la suite exacte duale de (\ref{segen}), on a
$A_{k,n}\subset (im\,f_k\otimes \Lambda^n)\cap (\Lambda^{k+1}\otimes
ker\,\theta_{k,n,1}).$ Posons à présent
$A'_{k,n}=A_{k,n}\cap(W_{(k+1,k)}\otimes\Lambda^n)$ : $A_{k,n}/A'_{k,n}$
s'injecte dans $(im\,f_{k}/W_{(k+1,k)})\otimes\Lambda^n$, qui n'a pas de
facteurs de composition du type mentionné dans (\ref{mri}) grâce à la proposition \ref{detfcl}.
Quant à $A'_{k,n}/W_{(k+1,k,n)}$, il s'injecte dans $\Lambda^{k+1}\otimes\big(\underset{t\geq 2}{\bigoplus}\Lambda^{k+t,n-t}\big)$, de sorte que le théorème \ref{fclambda} suffit à conclure.

L'assertion (\ref{mri7}) est analogue, en remarquant que la suite exacte
duale de (\ref{segen}) fournit $B_{k,n}\subset (im\,g_k\otimes\Lambda^n)\cap (\Lambda^{k+3}\otimes ker\,\theta_{k,n,1})$, parce que le morphisme 
$$\Lambda^{k+2}\otimes\Lambda^{k+1}\otimes\Lambda^n\xrightarrow{\theta\otimes\Lambda^n}\Lambda^{k+3}\otimes\Lambda^k\otimes\Lambda^n\xrightarrow{\Lambda^{k+3}\otimes\theta}\Lambda^{k+3}\otimes\Lambda^{k+1}\otimes\Lambda^{n-1}$$
est la somme des morphismes 
$$\Lambda^{k+2}\otimes\Lambda^{k+1}\otimes\Lambda^n\simeq\Lambda^{k+1,k+2,n}\xrightarrow{\Lambda^{k+1}\otimes\theta}\Lambda^{k+1,k+3,n-1}\simeq\Lambda^{k+3,k+1,n-1}$$
et
$$\Lambda^{k+2}\otimes\Lambda^{k+1}\otimes\Lambda^n\xrightarrow{\Lambda^{k+2}\otimes\theta}\Lambda^{k+2}\otimes\Lambda^{k+2}\otimes\Lambda^{n-1}\xrightarrow{\theta\otimes\Lambda^n}\Lambda^{k+3}\otimes\Lambda^{k+1}\otimes\Lambda^{n-1}.$$
\end{proof}

Nous terminons ce paragraphe en donnant une estimation de la division
par $\Lambda^1$ des foncteurs que nous avons introduits. 

On commence par observer que, comme les
foncteurs $L(2)$ et $\bar{D}(2)$ sont des quotients de
$\bar{I}^{\otimes 2}$, leur division par $\Lambda^1$ est nulle, de
sorte que $(L(2)\otimes\Lambda^n:\Lambda^1)\simeq
L(2)\otimes\Lambda^{n-1}$ et
$(\bar{D}(2)\otimes\Lambda^n:\Lambda^1)\simeq\bar{D}(2)\otimes\Lambda^{n-1}$.
La démonstration ci-dessous exploite sans cesse ces identifications et
d'autres analogues.

\begin{pr}\label{dluln}
L'image du morphisme $(L^2_n:\Lambda^1)\to
L(2)\otimes\Lambda^{n-1}$ (resp. $(D^2_n:\Lambda^1)\to
\bar{D}(2)\otimes\Lambda^{n-1}$) induit par l'inclusion est incluse dans $L^2_{n-1}$ (resp. $D^2_{n-1}$).
\end{pr}

\begin{proof} La suite exacte duale de (\ref{segen}) montre que la
  composée $L^2_n\oplus
  D^2_n\hookrightarrow\Lambda^2(\bar{I})\otimes\Lambda^n\xrightarrow{Dh_n}\bar{I}^{\otimes 2}\otimes\Lambda^{n-1}$ est nulle ; par division par $\Lambda^1$, on obtient que la somme directe des images des morphismes de l'énoncé est incluse dans le noyau de $\Lambda^2(\bar{I})\otimes\Lambda^{n-1}\xrightarrow{(Dh_n : \Lambda^1)}\bar{I}^{\otimes 2}\otimes\Lambda^{n-2}$. Il suffit donc de vérifier que $(Dh_n : \Lambda^1)=h_{n-1}$, ce qui provient de la commutation du diagramme
$$\xymatrix{\Lambda^2(\bar{I})\otimes\Lambda^n\ar[rr]^-{Dh_n}\ar@{^{(}->}[d]
  & & 
  \bar{I}^{\otimes 2}\otimes\Lambda^{n-1}\ar@{^{(}->}[d] \\
\Lambda^2(\bar{I})\otimes\Lambda^{n-1}\otimes\Lambda^1\ar[rr]^-{Dh_{n-1}\otimes\Lambda^1}
& & 
  \bar{I}^{\otimes 2}\otimes\Lambda^{n-2}\otimes\Lambda^1
}$$
qu'on obtient en dualisant le diagramme
$$\xymatrix{\Lambda^2(\bar{P})\otimes\Lambda^n
  & & 
  \bar{P}^{\otimes 2}\otimes\Lambda^{n-1}\ar[ll]_-{h_n} \\
\Lambda^2(\bar{P})\otimes\Lambda^{n-1}\otimes\Lambda^1\ar@{>>}[u]
& & 
  \bar{P}^{\otimes 2}\otimes\Lambda^{n-2}\otimes\Lambda^1\ar[ll]_-{h_{n-1}\otimes\Lambda^1}\ar@{>>}[u]
}$$
qui est commutatif par inspection.
\end{proof}

\begin{nota} Nous désignerons par $u^L_n : (L^2_n:\Lambda^1)\to
  L^2_{n-1}$ et $u^D_n : (D^2_n:\Lambda^1)\to D^2_{n-1}$ les
  morphismes procurés par la proposition précédente.
\end{nota}

\begin{rem} On peut montrer que ces flèches sont des
  isomorphismes. Leur surjectivité se déduit d'ailleurs aisément des
  considérations du paragraphe suivant.

Dans \cite{these}, nous établissons ce type de résultat dans un cadre plus général et conceptuel.
\end{rem}

\subsection{Démonstration du théorème principal}

La proposition suivante constitue la clef de voûte de notre approche
de la structure des foncteurs $I^{\otimes 2}\otimes\Lambda^n$. Les
théorèmes \ref{thp} et \ref{thld}, dont les énoncés seront précisés dans les
théorèmes \ref{parts} et \ref{thpp}, s'en déduit par des arguments formels, moyennant le
théorème \ref{lmgp} qui traite le cas $n=0$.

\begin{pr} \label{pcruc}
Soient $n\in\mathbb{N}^*$ et $X$ un sous-foncteur de $L^2_n$
(resp. $D^2_n$). On suppose que le morphisme $(X:\Lambda^1)\to (L^2_n:\Lambda^1)\xrightarrow{u^L_n}
L^2_{n-1}$ (resp. $(X:\Lambda^1)\to (D^2_n:\Lambda^1)\xrightarrow{u^D_n}D^2_{n-1}$) est surjectif. 

On a alors $X=L^2_n$ (resp. $X=D^2_n$).
\end{pr}

\begin{proof} On traite d'abord le cas de $L^2_n$. Pour tout entier naturel impair $i$, la partition 
  $(i+n+1,i+n,n)$ est alternée, donc Weyl-séparante par la proposition
  \ref{pfc}. Cela permet d'appliquer la proposition \ref{cr-fd1} au
  sous-objet $L^2_n$ de $L(2)\otimes\Lambda^n$, où l'on prend pour
  $\beta$ le morphisme $f_{i+n}$, $\alpha$ étant ensuite défini par la
  composition de la condition {\em 1} de l'énoncé de ladite
  proposition. La proposition \ref{lpcfc} montre que les deuxième et troisième conditions de la proposition
  \ref{cr-fd1} sont vérifiées \ref{lpcfc}. La dernière condition en
  est également satisfaite puisque toute partition de longueur au
  plus $2$ de $2(i+n)+1$ est supérieure à $(i+n+1,i+n)$, et que
  $ker\,f_{i+n}$ ne peut avoir de facteur de composition
  $S_{(i+n+1,i+n)}$, puisque sa source en possède un seul et que son
  image en a un par la proposition \ref{detfcl}. 

Par conséquent, la partition  $(i+n+1,i+n,n)$ est $\Lambda^1$-détectable dans
  $L^2_n$  relativement au morphisme $u^L_n$, lorsque l'entier naturel
  $i$ est impair. Ainsi, on a $(i+n+1,i+n,n)\vdash X$ pour $i$ impair,
  d'où $X=L^2_n$ par le corollaire \ref{cdetl}.

Le cas de $D^2_n$ se traite pareillement, en considérant la
partition alternée $(i+n+3,i+n,n)$ et le morphisme $g_{i+n}$ pour $i$
impair. La seule différence réside dans la satisfaction de la dernière
hypothèse de la proposition \ref{cr-fd1} : on doit utiliser que
$\bar{D}(2)$ n'a pas de facteur de composition $S_{(i+n+2,i+n+1)}$
(cf. proposition \ref{detfcl}), et que toute partition  de $2(i+n)+3$ de longueur au
plus $2$ et distincte de $(i+n+2,i+n+1)$ est supérieure à $(i+n+3,i+n)$.
\end{proof}

\begin{theo}\label{parts}
Pour tout $n\in\mathbb{N}$, $L^2_n$ et $D^2_n$ sont simples artiniens de type $2$.
\end{theo}

\begin{proof} Les foncteurs $L^2_n$ et $D^2_n$ ne sont pas artiniens
  de type $1$, car l'image par $Dg^L_n$ et $Dg^D_n$ des filtrations respectives de $L(2)$ et
  $\bar{D}(2)$ du théorème \ref{lmgp},
  tensorisées par $\Lambda^n(\Lambda^3)$, en fournit des filtrations
  infinies de quotients infinis. 

Montrons maintenant par récurrence sur $n$ que pour tout sous-objet strict $X$ de $L^2_n$ (resp. $D^2_n$) et tout foncteur fini $F$, $X\otimes F$ est artinien de type $1$. Pour $n=0$, cette assertion est incluse dans le théorème \ref{lmgp}.

Supposons maintenant $n>0$ et l'assertion démontrée pour $L^2_{n-1}$
(resp. $D^2_{n-1}$). Si $X$ est un sous-objet strict de $L^2_n$
(resp. $D^2_n$), la proposition \ref{pcruc} montre que l'image $A$ du morphisme
$f :(X:\Lambda^1)\to L(2)\otimes\Lambda^{n-1}$
(resp. $(X:\Lambda^1)\to \bar{D}(2)\otimes\Lambda^{n-1}$) induit par
l'injection de $X$ dans $L(2)\otimes\Lambda^n$
(resp. $\bar{D}(2)\otimes\Lambda^n$) est un sous-objet strict de $L^2_{n-1}$ (resp. $D^2_{n-1}$). 

Le diagramme commutatif 
$$\xymatrix{X\ar@{^{(}->}[r]\ar[d] & L(2)\otimes\Lambda^n
  \ar@{^{(}->}[d] \\
(X:\Lambda^1)\otimes\Lambda^1 \ar[r]^-{f\otimes\Lambda^1} & L(2)\otimes\Lambda^{n-1}\otimes\Lambda^1,
}$$
dans lequel la flèche verticale de gauche est l'unité de l'adjonction,
montre que $X$ s'injecte dans $A\otimes\Lambda^1$. Ce foncteur étant artinien de
type $1$, de même que son produit tensoriel par un foncteur fini, par hypothèse de récurrence, cela termine la démonstration.\end{proof}

\begin{theo}\label{thpp}
Pour tout $n\in\mathbb{N}$, $I^{\otimes 2}\otimes\Lambda^n$ est artinien de type $2$.
\end{theo}

\begin{proof}
On commence par remarquer que si $\Lambda^2(\bar{I})\otimes\Lambda^n$
est artinien de type $2$, il en est de même pour $\bar{I}^{\otimes
  2}\otimes\Lambda^n$, donc de $I^{\otimes 2}\otimes\Lambda^n$. Il suffit pour cela de considérer la filtration de $\bar{I}^{\otimes 2}$ de sous-quotients $\Lambda^2(\bar{I})$, $\bar{I}$ et $\Lambda^2(\bar{I})$, sachant que $\bar{I}\otimes \Lambda^n$ est artinien de type $1$ (cf. \cite{Piriou}). 

On procède ensuite par récurrence sur $n$, le cas $n=0$ étant rappelé
dans la proposition \ref{lmgp}. Supposons donc $n>0$ et
$\bar{I}^{\otimes 2}\otimes\Lambda^{n-1}$ artinien de type $2$. La
suite exacte duale de (\ref{segen})  montre que
$(\Lambda^2(\bar{I})\otimes\Lambda^n)/(L^2_n\oplus\bar{D}^2_n)$
s'injecte dans $\bar{I}^{\otimes 2}\otimes\Lambda^{n-1}$, ce quotient
est donc artinien de type $2$. Le théorème précédent implique donc le résultat.
\end{proof} 

\begin{rem} On peut retrouver les résultats de Piriou (\cite{Piriou})
  relatifs aux foncteurs $\bar{I}\otimes\Lambda^n$ par la même
  méthode. L'article \cite{Piriou} repose également sur l'étude de
  facteurs de composition idoines, mais procède de manière beaucoup
  plus explicite, à l'aide de calculs de groupes d'extensions. Powell,
  qui a généralisé dans \cite{Po3} les résultats de
  Piriou au cas du produit tensoriel entre $\bar{I}$ et un foncteur
  fini, mène des raisonnements sur des facteurs de composition
à l'aide de quotients du foncteur différence, dont le
maniement est cependant différent de celui de $(-:\Lambda^1)$, de
sorte que sa stratégie globale, tout en présentant des similitudes
avec celle du présent article, en diverge conceptuellement.
\end{rem}

%\dem on raisonne par récurrence sur $n$, le cas $n=0$ étant rappelé dans la proposition \ref{lmgp}. Supposons donc $n>0$ et $\bar{I}^{\otimes 2}\otimes\Lambda^{n-1}$ artinien de type $2$. La suite exacte duale de (\ref{segen}) montre (en utilisant aussi le point \ref{ptdg2} de la proposition \ref{ptdg}) que $(\Lambda^2(\bar{I})\otimes\Lambda^n)/(L^2_n\oplus\bar{D}^2_n)$ s'injecte dans $\bar{I}^{\otimes 2}\otimes\Lambda^{n-1}$, ce quotient est donc artinien de type $2$. La proposition précédente implique donc que $\Lambda^2(\bar{I})\otimes\Lambda^n$ est artinien de type $2$. Comme $\bar{I}^{\otimes 2}\otimes\Lambda^n$ admet une filtration dont les sous-quotients sont $\Lambda^2(\bar{I})\otimes\Lambda^n$, $\bar{I}\otimes\Lambda^n$, qui est artinien de type $1$, et $\Lambda^2(\bar{I})\otimes\Lambda^n$, le corollaire est démontré.

\paragraph*{Remerciements} L'auteur tient à témoigner sa gratitude
envers Lionel Schwartz pour ses nombreuses discussions sur les modules
instables et la catégorie $\F$, ainsi qu'à Geoffrey Powell, tant pour ses
remarques mathématiques que pour ses commentaires qui ont grandement contribué
à améliorer la présentation de cet article. Il remercie également Christine Vespa pour ses
encouragements et ses conseils.

\nocite{*}
\bibliographystyle{smfalpha}
\bibliography{artb}

\begin{flushleft}
\small{Aurélien DJAMENT\\
LAGA, Institut Galil\'ee\\
universit\'e Paris 13\\
99 avenue J.-B. Clément\\
93430 VILLETANEUSE (FRANCE)\\

\smallskip

djament@math.univ-paris13.fr}
\end{flushleft}
\end{document}